\documentclass[11pt]{article}
\usepackage{amsmath,amsfonts}
\usepackage{xcolor}
\usepackage{amsthm}
\usepackage{dsfont}			
\usepackage{hyperref}
\usepackage{multirow}
\usepackage{graphicx}
\usepackage{subcaption}
\usepackage{epsfig}
\def\cb{{\bf c}}
\def\tb{{\bf t}}

\def\l2{{l_2}}

\newcommand*{\Ztuk}{{}_{n}\mathring{\mathbb{Z}}_{u_k}}
\newcommand*{\thetab}{{\boldsymbol{\theta}}}
\newcommand*{\etab}{{\boldsymbol{\eta}}}

\newcommand*{\lb}{{\boldsymbol{\ell}}}
\newcommand*{\E}[1]{\operatorname{E}_{#1}}

\newcommand*{\deq}{{\, \stackrel{\text{def}}{=} \; }}
\newtheorem{theorem}{Theorem}

\hypersetup{
	pdfborder = {0 0 0},
	urlbordercolor = {0 0 0},
	colorlinks = {false},
}

\newcommand\smallO{
	{{\scriptscriptstyle\mathcal{O}}}
}
\def\underacc #1/#2{\mathchoice{\uacc\textstyle{#1}{#2}}{\uacc\textstyle{#1}{#2}}
	{\uacc\scriptstyle{#1}{#2}}{\uacc\scriptscriptstyle{#1}{#2}}}
\def\uacc#1#2#3{\mathop{#2{}}\limits_{#1#3{}}}

\newcommand*{\ind}[1]{\mathds{1} \left( {#1} \right)}

\begin{document}
	\title{Adaptive exact recovery in sparse nonparametric models}
	\author{Natalia Stepanova$^a$
		and Marie Turcicova$^{b0}$\\
		\\
		\small {\textit{$^a$School of Mathematics and Statistics, Carleton University,}}\\
		\small {\textit{1125 Colonel By Drive, Ottawa, ON, K1S 5B6, Canada} } \\
		\small{\textit{$^b$Institute of Computer Science, Czech Academy of Sciences, }} \\
		\small{\textit{Pod Vodárenskou věží 271/2, Prague, 182 00,  Czech Republic}}}
	\date{}
	
	\footnotetext{Corresponding author. Email address: turcicova@cs.cas.cz}
	\maketitle

\begin{abstract}
We observe an unknown  function of $d$ variables $f(\tb)$, $\tb\in[0,1]^d$, in the Gaussian white noise model of intensity $\varepsilon>0$.
We assume that the function $f$ is regular and that it is a sum of $k$-variate functions, where $k$ varies from $1$ to $s$ ($1\leq s\leq d$).
These functions are unknown to us and only a few of them are nonzero.
In this article, we address the problem of identifying the nonzero components of $f$  in the case when $d=d_\varepsilon\to \infty$ as $\varepsilon\to 0$ and $s$ is either fixed or
$s=s_\varepsilon\to \infty$, $s=o(d)$ as $\varepsilon\to \infty$. This may be viewed as a variable selection problem.
We derive the conditions when exact variable selection in the model at hand is possible and provide a selection procedure
that achieves this type of selection. The procedure is adaptive to a degree of model sparsity
described by the sparsity parameter $\beta\in(0,1)$.
We also derive conditions that  make the exact variable selection impossible. Our results augment previous work in this area.
\textit{2020MSC}: Primary: 62G08, Secondary: 62H12, 62G20
\end{abstract}
	
	\section{Introduction} \label{sec:problem_statement}
Analysis of high-dimensional data has recently become more and more important since large datasets have become increasingly available in every field of research -- ranging from genomic sequencing and brain imaging to climate monitoring and social sciences.
However, such analysis is often associated with specific phenomena that go beyond classical statistical theory.
The main difficulty is that the dimension can be very large. Speaking in mathematical terms, it is natural to
assume that the dimension tends to infinity.

In this work, we address the problem of adaptive recovery of the sparsity pattern
of a  multivariate signal observed in the Gaussian white noise,
and augment the results obtained earlier in Ingster and Stepanova (2014) and Stepanova and Turcicova (2025).
Specifically, we assume that an unknown signal $f$ of $d$ variables is observed in the Gaussian white noise model
\begin{gather}\label{model0}
dX_\varepsilon(\tb)=f(\tb)d\tb+\varepsilon d W(\tb),\quad \tb\in[0,1]^d,
\end{gather}
where $dW$ is a $d$-parameter Gaussian white noise and $\varepsilon>0$ is the noise intensity.
 The signal 	$f$ belongs to a subspace of $L_2([0,1]^d)=L_2^d$ with inner product  $(\cdot,\cdot)_2$ and norm $\|\cdot\|_2$
  that consists
 of regular enough functions, and we assume that $d=d_\varepsilon\to \infty$ as $\varepsilon\to 0$.
 Consider an operator $\mathbb{W}:L_2^d\to {\cal G}_0$ taking values in
 the set ${\cal G}_0$ of centered Gaussian random variables
 such that if $\xi_0=\mathbb{W}(g_1)$ and $\eta_0=\mathbb{W}(g_2)$, where $g_1,g_2\in L_2^d$,
then ${\rm cov}(\xi_0,\eta_0)=(g_1,g_2)_{2}$. The $d$-parameter Gaussian white noise $dW$ in model~(\ref{model0}) is defined through the operator $\mathbb{W}$ by
  $$\mathbb{W}(g)=\int_{[0,1]^d} g(\tb) dW(\tb)\sim N(0,\|g\|_2^2),\quad g\in L_2^d.$$
  In particular, if $\{g_\lb\}_{\lb\in{\cal L}}$ is an  orthonormal basis of $L_2^d$, then $\mathbb{W}(g_\lb)\sim N(0,1)$ for $\lb\in{\cal L}$
 and, for any finite set $\{g_{\lb}\}$ of the basis functions, the family $\{\mathbb{W}(g_{\lb})\}$ forms a multivariate standard normal vector. Thus,
  the centered Gaussian measure on $L_2^d$ determined by $\mathbb{W}$ has a diagonal covariance operator (i.e., the identity operator).
  Furthermore, let $\mathbb{X}_\varepsilon:L_2^d\to {\cal G}$ be an operator taking values in the set
  ${\cal G}$  of Gaussian random variables such that if $\xi=\mathbb{X}_{\varepsilon}(g_1)$ and $\eta=\mathbb{X}_{\varepsilon}(g_2)$, where $g_1,g_2\in L_2^d$,
then ${\rm E}(\xi)=(f,g_1)_2$, ${\rm E}(\eta)=(f,g_2)_2$, and ${\rm cov}(\xi,\eta)=\varepsilon^2(g_1,g_2)_2$.
 By ``observing the trajectory (\ref{model0})'' we mean observing
 a realization of the Gaussian field $X_\varepsilon(\tb)$, $\tb\in[0,1]^d$,  defined
 through the operator $\mathbb{X}_\varepsilon$
 by
 \begin{gather*}
 \mathbb{X}_\varepsilon(g)=\int_{[0,1]^d}g(\tb)dX_\varepsilon(\tb)\sim N\left((f,g)_2,\varepsilon^2\|g\|_2^2  \right),\quad  g\in L_2^d.
 \end{gather*}
  In terms of the operators $\mathbb{W}$ and $\mathbb{X}_\varepsilon$, the stochastic differential equation (\ref{model0}) can be expressed as
\begin{equation}\label{model1}
		\mathbb{X}_{\varepsilon}=f+\varepsilon \mathbb{W},
	\end{equation}
and ``observing the trajectory (\ref{model1})'' means that we observe all normal $N\left((f,g)_2,\varepsilon^2\|g\|_2^2  \right)$ random variables
when $g$ runs through $L_2^d$. For any $f\in L_2^d$, the ``observation'' $\mathbb{X}_\varepsilon$ in model (\ref{model1})
defines the Gaussian measure ${\rm P}_{\varepsilon,f}$ on the Hilbert space $L_2^d$ with mean function $f$ and covariance operator
$\varepsilon^2I$, where $I$ is the identity operator (for references, see Gin\'{e} and Nickl (2016), Ibragimov and Khasminskii (1997),
Skorohod (1998)).
In addition to regularity constraints, we assume that $f$ has a sparse structure and consider the problem of recovering exactly the
sparsity pattern of $f$ by using the asymptotically minimax approach.

\subsection{Sparsity conditions}
In high-dimensional inference problems like the one under study,
the quality of statistical procedures deteriorates very quickly as the dimension increases.
One way to avoid the curse of dimensionality stemming from high-dimensional settings is to reduce the ``working dimension'' of the problem.
In the present setup, this can be done by employing functional ANOVA-type decompositions for~$f$.
One common functional ANOVA model (see, for example, Lin (2000), Ingster and Suslina (2015)) assumes that
the function of $d$ variables is a sum of functions of one variable (main effects),
functions of two variables (two-way interactions), and so on. That is, the function $f$ of $d$ variables is decomposed as
 \begin{equation*}
 	\begin{split}
	f(\tb)&=\sum_{u\subseteq\{1,\ldots,d\}} f_u(\tb_u),\quad \tb=(t_1,\ldots,t_d)\in[0,1]^d, \quad \tb_{u}=(t_j)_{j \in u},
	\label{fun_anova}
	 	\end{split}
\end{equation*}
where the sum is taken over all subsets $u \subseteq \{1,\ldots,d\}$. Owen (1998) noted that the functional ANOVA decomposition is completely analogous to the one used in experimental statistics because the total variance $\sigma^2 = \int_{[0,1]^d} (f(\tb) - {\mathds I})^2 d\tb$, where ${\mathds I} = \int_{[0,1]^d} f(\tb) d\tb$, can be written as $\sigma^2 = \sum_{u\subseteq \{1,\ldots,d\}} \sigma^2_u$, where $\sigma_u^2 = \int_{[0,1]^k} f_u^2(\tb_u) d \tb_u$ if $u \neq \emptyset$ and $\sigma^2_{\emptyset} = 0$.
In Ingster and Suslina (2015), this decomposition is used to set up nonparametric alternatives to the hypothesis of no signal in a signal detection problem. In~Owen (1998) and Lin (2000), the functional ANOVA decomposition is used to split the problem of integrating and estimating a multivariate function $f$ into low-dimensional tasks.

 It is assumed that, in the above functional ANOVA decomposition, $f_u={\rm const}$ for $u=\emptyset$  and, in order to guarantee uniqueness, that
\begin{gather*}
	\int_{0}^1 f_u(\tb_u)\,dt_j=0,\quad \mbox{for}\; j\in u,
\end{gather*}
that is, the terms $f_u$ are mutually orthogonal in $L_2^d$.
Each function $f_u(\tb_u)$ depends only on variables in $\tb_u$ and describes the ``interaction'' between these variables.
Denote by $\#(u)$ the number of elements of a set $u$ and let $\tb_{-u} = (t_j)_{j \notin u}$.
One example of the function $f_u(\tb_u)$ that  satisfies the above orthogonality condition is
(see Owen (1998) for details)
$$ f_u (\tb_u) = \int_{[0,1]^{d-\#(u)}} \left( f(\tb) - \sum_{v \subset u\atop v\neq u} f_v(\tb_v)\right) d\tb_{-u}.$$
Another example is the function $f_u(\tb_u)=\prod_{j\in u}f_j(t_j)$, where each  $f_j$ integrates to zero.

In high-dimensional inference problems, alternatively, or in addition to the general functional ANOVA model,
we may choose some $s$ ($1\leq s\leq d$)
and consider the expansion
 \begin{equation*}
		f(\tb)=\sum_{ \substack{u\subseteq\{1,\ldots,d\},\, 1\leq \# (u) \leq s}} f_u(\tb_u),\quad \tb\in[0,1]^d,
\label{fun_sub_anova}
\end{equation*}
where, if $s$ is small relative to $d$, $f(\tb)$ is a sum of  functions of a small number of variables.
The model of interest in this work is obtained from the above orthogonal expansion
by assuming additionally that $f(\tb)$ has a sparse structure.
To ensure that the recovery of the $d$-variate signal $f$ in model (\ref{model1}) is feasible,
the sparsity and regularity
constraints on $f$ are required.

For $1\leq k\leq s$, where $s$ is as above, define ${\cal U}_{k,d}=\{u_k:u_k\subseteq\{1,\ldots,d\}, $ $ \#(u_k)=k\}$ and observe that $\# (\mathcal{U}_{k,d})={d\choose k}.$
If $u_k=\{j_1,\ldots,j_k\}\in {\cal U}_{k,d}$, $1 \leq j_1 < \ldots < j_k \leq d$, denote as before
$\tb_{u_k}=(t_{j_1},\ldots,t_{j_k})\in[0,1]^k$, and assume that the signal $f$ in model (\ref{model1}) has the form
 \begin{equation}
f(\tb)=\sum_{k=1}^{s} \sum_{u_k\in{\cal U}_{k,d}}\eta_{u_k} f_{u_k}(\tb_{u_k}),\quad \tb\in [0,1]^d,
	\label{fun_sub_anova22}
\end{equation}
where the components $f_{u_k}$, $ u_k\in {\cal U}_{k,d}, $ $1\leq k\leq s,$ satisfy the orthogonality condition
\begin{gather}\label{orthcon2}
	\int_{0}^1 f_{u_k}(\tb_{u_k})\,dt_j=0,\quad \mbox{for}\; j\in u_k,
\end{gather}
the nonrandom quantities $\eta_{u_k}\in\{0,1\}$, called the \textit{index variables}, mark the component $f_{u_k}$ as ``active''
when $\eta_{u_k}=1$ and as ``inactive'' when $\eta_{u_k}=0$, and as  $\varepsilon\to 0 $
 \begin{gather}\label{SC1}
\sum_{u_k\in{\cal U}_{k,d}}\eta_{u_k}=  {d\choose k}^{1-\beta}(1+o(1)),\quad 1\leq k\leq s,
 \end{gather}
where $\beta\in(0,1)$ is an unknown model parameter called the \textit{sparsity index}. Condition (\ref{SC1}) is the sparsity condition that determines
the sparsity structure of $f$ in (\ref{fun_sub_anova22}). This condition implies that if $1\leq k<l\leq s$,
the number of $k$-variate terms
on the right-hand side of (\ref{fun_sub_anova22}) that are active is less than the number of $l$-variate terms that are active.

In a simpler setup, when instead of decomposition (\ref{fun_sub_anova22}) in model (\ref{model1})--(\ref{SC1}), the signal $f(\tb)$ is assumed to be a sum of functions, each depending on $k$ variables ($1\leq k\leq d$), of the
form
\begin{equation}\label{f}
	f(\tb)=\sum_{u_k\in{\cal U}_{k,d}}\eta_{u_k} f_{u_k}(\tb_{u_k}),\quad \tb\in[0,1]^d,
\end{equation}
where  $\etab_k=(\eta_{u_k}, u_k\in {\cal U}_{k,d})$ is such that $\sum_{u_k\in{\cal U}_{k,d}}\eta_{u_k}=  {d\choose k}^{1-\beta}(1+o(1))$
as $\varepsilon\to 0$,
the problem of identifying exactly the active components of $f(\tb)$ was studied in Ingster and Stepanova (2014) for $k=1$ and
in Stepanova and Turcicova (2025) for $1\leq k\leq d$. In these articles, conditions under which exact identification of the nonzero components
of $f(\tb)$ as in (\ref{f}) is possible and impossible were established,
and adaptive selection procedures that, under some model assumptions, identify exactly all nonzero $k$-variate components $f_{u_k}$
were proposed.

In this work, we consider the two cases: (i) when $s$ is fixed and (ii) when $s=s_\varepsilon\to \infty$, $s=o(d)$ as $\varepsilon\to 0$. In both cases,
the sparsity condition (\ref{SC1}) implies
\begin{equation*}
	\sum_{k=1}^{s} \sum_{u_k \in \mathcal{U}_{k,d}} \eta_{u_k}= \sum_{k=1}^s {d\choose k}^{1-\beta}(1+o(1))= {d\choose s}^{1-\beta}(1+o(1)),
\quad \varepsilon\to 0,
\end{equation*}
that is, only ${d\choose s}^{1-\beta}(1+o(1))=o({d\choose s})$ orthogonal components $f_{u_k}$ of $f$ in (\ref{fun_sub_anova22}) are active
and the remaining
components are inactive, which implies that the function $f$ is \textit{sparse}.
The values of $\beta$ that are close to 1 make the signal $f$ in (\ref{fun_sub_anova22}) \textit{highly sparse}, whereas the
values of $\beta$ that are close to 0 make it \textit{dense}. Condition (\ref{SC1}), in which $1\leq s\leq d$, may be viewed as a natural
extension of the sparsity condition
 \begin{gather*}
 \sum_{u_k\in{\cal U}_{k,d}}\eta_{u_k}=  {d\choose k}^{1-\beta}(1+o(1))\quad \varepsilon\to 0,\quad \mbox{for some}\;1\leq k\leq d,
 \end{gather*}
 employed in Stepanova and Turcicova (2025) to recover the sparsity pattern of function $f(\tb)$ as in (\ref{f}) for a chosen $k$.

Define the sets
${H}_{\beta,d}^k={H}_{\beta,d}^k(\varepsilon)$, $1\leq k\leq s$, and
$\mathcal{H}^{s}_{\beta,d} =\mathcal{H}^{s}_{\beta,d}(\varepsilon)$ as follows:
\begin{align*}
{H}_{\beta,d}^k &= \left\{ \boldsymbol{\eta}_k=( \eta_{u_k}, u_k \in \mathcal{U}_{k,d}) : \eta_{u_k} \in \{0,1\},\, u_k\in {\cal U}_{k,d},\; \sum\nolimits_{u_k\in{\cal U}_{k,d}}\eta_{u_k}= {d\choose k}^{1-\beta}(1+o(1)) \right\},\\
{\cal H}_{\beta,d}^{s}&=\{\etab=(\etab_1,\ldots,\etab_s): \etab_k\in { H}^k_{\beta,d},\,1\leq k\leq s\}.
\end{align*}
Based on the ``observation'' $\mathbb{X}_{\varepsilon}$ in model (\ref{model1})--(\ref{SC1}), we wish to identify,
with high degree of accuracy, the sparsity pattern of $f$, that is,
we wish to construct a good  estimator $\boldsymbol{\hat{\eta}}=(\boldsymbol{\hat{\eta}}_1,\ldots,\boldsymbol{\hat{\eta}}_s)$ of
$\etab=(\etab_1,\ldots,\etab_s)\in \mathcal{H}^{s}_{\beta,d}$ that would tell us which terms $f_{u_k}$ in the sparse decomposition (\ref{fun_sub_anova22}) are active. This may be viewed as a variable selection
problem, and $\hat{\etab}$ may be named a \textit{selector}.

\subsection{Regularity conditions}\label{RC}
Attempting to provide an asymptotically minimax solution to the  problem of recovering the sparsity pattern of a function~$f$ in model  (\ref{model1})--(\ref{SC1}),
we have to assume that the set of signals $f$  is not ``too large".
In this article, we will be interested  in periodic Sobolev classes described by means of Fourier coefficients.
Such classes are quite common in the literature on nonparametric estimation, signal detection, and variable selection.

Namely, following the construction in Ingster and Suslina (2015),
for ${u_k}\in{\cal U}_{k,d}$, $1\leq k\leq d,$ consider the set
\begin{eqnarray*}
		\mathring{\mathbb{Z}}_{u_k}&=&\{\lb=(l_1,\ldots,l_d)\in \mathbb{Z}^d: l_j=0\;\mbox{for}\;j\notin u_k \mbox{ and }  l_j\neq 0\;\mbox{for}\;j\in u_k \}.
\end{eqnarray*}
We also set $\mathring{\mathbb{Z}}_\emptyset=\underbrace{(0,\ldots,0)}_d$, $\mathring{\mathbb{Z}}=\mathbb{Z}\setminus\{0\},$
${\mathring{\mathbb{Z}}}^k=\underbrace{\mathring{\mathbb{Z}}\times\ldots\times \mathring{\mathbb{Z}}}_k,$ and note that (for $k=0$ we set ${\cal U}_{k,d}=\emptyset$)
$$\mathbb{Z}^d=\left( \mathring{\mathbb{Z}}\cup\{0\} \right)^d= \bigcup_{u\subseteq\{1,\ldots,d\}}\mathring{\mathbb{Z}}_{u}=
\bigcup_{0\leq k\leq d}\bigcup_{u_k\in {\cal U}_{k,d}}\mathring{\mathbb{Z}}_{u_k}.$$
Consider the Fourier basis $\{\phi_{\lb}(\tb)\}_{\lb\in	\mathbb{Z}^d}$ of $L_2^d$ defined as follows:
\begin{equation}
	\begin{split}
		\phi_{\lb}(\tb)&=\prod_{j=1}^d \phi_{l_j}(t_{j}),\quad \lb=(l_1,\ldots,l_d)\in {\mathbb{Z}}^d,
		\\
		\phi_0(t)=1,\quad \phi_l(t)&=\sqrt{2}\cos(2\pi l t),\quad \phi_{-l}(t)=\sqrt{2}\sin(2\pi l t),\quad l>0.
		\label{OB}
	\end{split}
\end{equation}
Observe that
$\phi_{\lb}(\tb)=\phi_{\lb}(\tb_{u_k})$ for $\lb\in\mathring{\mathbb{Z}}_{u_k}$  (for $u=\emptyset$ we set $\phi_{\lb}(\tb_u)=1$) and
\begin{gather}\label{decomp0}
\{\phi_{\lb}(\tb)\}_{\lb\in\mathbb{Z}^d}=
\bigcup_{u\subseteq\{1,\ldots,d\}}\{\phi_{\lb}(\tb_{u})\}_{\lb\in	\mathring{\mathbb{Z}}_{u}}=
\bigcup_{0\leq k\leq d}\bigcup_{u_k\in {\cal U}_{k,d}}\{\phi_{\lb}(\tb_{u_k})\}_{\lb\in	\mathring{\mathbb{Z}}_{u_k}}.
\end{gather}
Next, let $\theta_{\lb}(u_k)=(f_{u_k},\phi_{\lb})_{L_2^d}$ be the $\lb$th Fourier coefficient of $f_{u_k}$ for $\lb \in 	 \mathring{\mathbb{Z}}_{u_k}$, $u_k\in{\cal U}_{k,d}$, $1\leq k\leq d$.
	Then, for $u_k=\{j_1,\ldots,j_k\}\in {\cal U}_{k,d}$, where $1 \leq j_1 < \ldots < j_k \leq d$, $1\leq k\leq s$, the $k$-variate component $f_{u_k}$ on the right-hand side of (\ref{fun_sub_anova22}) can be expressed as
	\begin{equation*}
		f_{u_k}(\tb_{u_k})=\sum_{\lb\in \mathring{\mathbb{Z}}_{u_k}}\theta_{\lb}(u_k) \phi_{\lb}(\tb_{u_k}),
	\end{equation*}
	and the entire function $f$ in (\ref{fun_sub_anova22}) takes the form
\begin{gather}\label{f-decomp}
	f(\tb)=\sum_{k=1}^{s} \sum_{u_k\in {\cal U}_{k,d}}\eta_{u_k} \sum_{\lb\in \mathring{\mathbb{Z}}_{u_k}}\theta_{\lb}(u_k) \phi_{\lb}(\tb_{u_k}).
\end{gather}
Note that only those Fourier coefficients of $f$ that correspond to the orthogonal components $f_{u_k}$ are nonzero and that
$\|f_{u_k}\|^2_{2}=(f_{u_k},f_{u_k})_{L_2^d}=\sum_{\lb\in \mathring{\mathbb{Z}}_{u_k}} \theta_{\lb}^2 (u_k)$.

  For $u_k=\{j_1,\ldots,j_k\}\in{\cal U}_{k,d},$ where $1\leq j_1<\ldots<j_k\leq d$, $1\leq k\leq s$, we assume that $f_{u_k}$ belongs to the Sobolev class of $k$-variate functions with integer smoothness parameter $\sigma\geq 1$ for which
 the semi-norm $\|\cdot\|_{\sigma,2}$ is defined by
 \begin{gather}\label{semi-norm}
 \|f_{u_k}\|_{\sigma,2}^2=\sum_{i_1=1}^k\ldots\sum_{i_{\sigma}=1}^k \left\|\frac{\partial^{\sigma} f_{u_k}  }{\partial t_{j_{i_1}}\ldots\partial t_{j_{i_{\sigma}}  }}   \right\|^2_2.
 \end{gather}
 Under the periodic constraint, we can define the semi-norm $\|\cdot\|_{\sigma,2}$ for the general case $\sigma>0$
 in terms of the Fourier coefficients $\theta_{\lb}(u_k)$, $\lb \in \mathring{\mathbb{Z}}_{u_k}$.
  For this, assume that $f_{u_k}(\tb_{u_k})$ admits 1-periodic $[\sigma]$-smooth extension in each argument to $\mathbb{R}^k$,
 i.e., for all derivatives $f_{u_k}^{(n)}$ of integer order $0\leq n\leq [\sigma]$, where $f_{u_k}^{(0)}=f_{u_k}$, one has
 \begin{gather*}
 	f_{u_k}^{(n)}(t_{j_1},\ldots,t_{j_{i-1}},0,t_{j_{i+1}},\ldots,t_{j_k})=f_{u_k}^{(n)}(t_{j_1},\ldots,t_{j_{i-1}},1,t_{j_{i+1}},\ldots,t_{j_k}),\quad 2\leq i\leq k-1,
 \end{gather*}
 with obvious extension for $i=1,k$. Then, the expression in (\ref{semi-norm}) corresponds to
 \begin{gather}
 	\|f_{u_k }\|_{\sigma,2}^2=\sum_{\lb\in \mathring{\mathbb{Z}}_{u_k} } \theta_{\lb}^2(u_k) c_{\lb}^2,\quad c_{\lb}^2=\left(\sum_{j=1}^d (2\pi l_{j})^2  \right)^{\sigma}
 	=\left(\sum_{i=1}^k (2\pi l_{j_i})^2  \right)^{\sigma}.
 	\label{def:c_lb}
 \end{gather}
Next, denote by ${\cal F}_{\cb_{u_k}}$ the Sobolev ball of radius 1 with coefficients ${\cb}_{u_k}=({c}_{\lb})_{\lb\in\mathring{\mathbb{Z}}_{u_k} }$, that is,
 \begin{gather*}
 	{\cal F}_{\cb_{u_k}}=\left \{  f_{u_k}(\tb_{u_k})= \sum_{\lb\in\mathring{\mathbb{Z}}_{u_k} } \theta_{\lb}(u_k)\phi_{\lb}(\tb_{u_k}), \; \tb_{u_k} \in [0,1]^k:
 	\sum_{\lb\in \mathring{\mathbb{Z}}_{u_k} } \theta_{\lb}^2(u_k) c_{\lb}^2 \leq 1 \right \}
 \end{gather*}
and assume that every component $f_{u_k}$ of $f$ in (\ref{fun_sub_anova22}) belongs to this Sobolev ball, that is,
 \begin{gather}\label{Fu}
 f_{u_k} \in {\cal F}_{\cb_{u_k}},\quad u_k \in {\cal U}_{k,d},\;1 \leq k \leq s.
 \end{gather}
 Thus, the model under study is specified by equations (\ref{model1})--(\ref{SC1}) and (\ref{Fu}).

The collections $\cb_{u_k}=(c_{\lb})_{\lb\in\mathring{\mathbb{Z}}_{u_k}}$, $u_k\in {\cal U}_{k,d}$,
 are invariant with respect to a particular choice of the elements of $u_k\in{\cal U}_{k,d}$ for all $1\leq k\leq d$. Hence, if $u_k=\{j_1,j_2,\ldots,j_k\}$, $1\leq j_1<\ldots<j_k\leq d$, and $\lb\in \mathring{\mathbb{Z}}_{u_k}$, we can write $c_{\lb}=c_{l_{j_1},\ldots,l_{j_k},0,\ldots,0}$.
Thinking of a transformation from $\mathring{\mathbb{Z}}^k$ to $\mathring{\mathbb{Z}}_{u_k}$ that maps
$\lb^{\prime}=(l^{\prime}_{1},\ldots,l^{\prime}_{k})\in \mathring{\mathbb{Z}}^k$ to
$\lb=(l_{j_1},\ldots,l_{j_k},0,\ldots,0)\in \mathring{\mathbb{Z}}_{u_k}$ according to the rule
$l_j = 0$ for $ j \notin u_k$ and $l_{j_p}=l^{\prime}_p$ for $p=1,\ldots,k,$
we may, if needed, represent each collection  $\cb_{u_k}=(c_{\lb})_{\lb\in\mathring{\mathbb{Z}}_{u_k}}$
as $\cb_k=(c_{\lb})_{\lb\in\mathring{\mathbb{Z}}^k}$. Clearly, the Sobolev balls ${\cal F}_{\cb_{u_k}}$ are isomorphic for all $u_k$
of cardinality $k$ ($1\leq k\leq d$).

 \subsection{Problem statement}
 Based on the ``observation'' $\mathbb{X}_\varepsilon$ in model~(\ref{model1})--(\ref{SC1}) and (\ref{Fu}),
 we wish to identify, with a high degree of accuracy, the active components $f_{u_k}$ of $f$.
This leads to the problem of obtaining an  estimator $\boldsymbol{\hat{\eta}}=(\boldsymbol{\hat{\eta}}_1,\ldots, \boldsymbol{\hat{\eta}}_s)$
for $\etab=(\etab_1,\ldots,\etab_s)\in \mathcal{H}^s_{\beta,d}$,
which may be viewed as a variable selection problem, since $\boldsymbol{\hat{\eta}}$ selects important variables $t_j, \, j=1,\ldots,d$, and their relevant interactions
in the ANOVA-type expansion (\ref{fun_sub_anova22}). The estimator $\boldsymbol{\hat{\eta}}$ is named a \textit{selector}
since it serves to select the active components of~$f$ in decomposition (\ref{fun_sub_anova22}).

Identifying the active components of $f$ in model (\ref{model1})--(\ref{SC1}) and (\ref{Fu}) is feasible  when
the components~$f_{u_k}$ of $f$ are not ``too small'', i.e., they are separated from zero in an appropriate way.
Therefore, for a given $u_k\in {\cal U}_{k,d}$, $1\leq k\leq s$ and $r >0$, we define the set
$$\mathring{\cal F}_{\cb_{u_k}}(r)=\{ f_{u_k} \in \mathcal{F}_{\cb_{u_k}}: \| f_{u_k}  \|_2\geq r\} $$
   and consider testing
     \begin{gather}
 	\mathbb{H}_{0,u_k} :  f_{u_k}=0 \quad\mbox{vs.}\quad \mathbb{H}^{\varepsilon}_{1,u_k}:  f_{u_k}\in \mathring{\cal F}_{\cb_{u_k}}(r_{\varepsilon,k}),\quad  u_k\in {\cal U}_{k,d},
 \quad 1\leq k\leq s,
 	\label{hypotheses0}
 \end{gather}
  where  $r_{\varepsilon,k}>0$, $1\leq k\leq s$, and $\max_{1\leq k\leq s}r_{\varepsilon,k}\to 0$ as $\varepsilon\to 0$.

 To ensure that exact recovery of the
sparsity pattern of $f$ (in the sense defined precisely at the end of this section)  is possible, we require the components $f_{u_k}$, $u_k\in {\cal U}_{k,d}$, $1\leq k\leq s$, to be
 \textit{selectable}, that is, we require the family (indexed by $\varepsilon$) of collections $r_\varepsilon=\{r_{\varepsilon,k}, 1\leq k\leq s\}$ in the hypothesis testing problem (\ref{hypotheses0}) to be above
  the so-called \textit{sharp selection boundary}.
Establishing the sharp selection boundary that makes the exact recovery of $f$
in model (\ref{model1})--(\ref{SC1}) and (\ref{Fu}) feasible is one of the goals in this article.
For this purpose, the hypothesis testing problem (\ref{hypotheses0}) will be employed as an auxiliary problem.

Since the Sobolev balls $ {\cal F}_{\cb_{u_k}}$ are isomorphic for all $u_k$ of cardinality $k$,
 $1\leq k\leq s$,
 the problem of testing  $\sum_{k=1}^s {d\choose k}$ pairs of the null and alternative hypotheses  in (\ref{hypotheses0}) reduces to that of testing
 \begin{gather}\label{HT1}
 	\mathbb{H}_{0,u_k} :  f_{u_k}=0 \quad \mbox{vs.}\quad \mathbb{H}^{\varepsilon}_{1,u_k}:  f_{u_k}\in \mathring{\cal F}_{\cb_{u_k}}(r_{\varepsilon,k}),\quad
 \mbox{for some}\;\; u_k\in{\cal U}_{k,d},\;\;  1\leq k\leq s,
 \end{gather}
 and the number of tests one needs to carry out decreases from $\sum_{k=1}^s {d\choose k}$ to $s$.
  For every chosen $u_k$, $1\leq k\leq s$, the hypotheses $\mathbb{H}_{0,u_k}$ and $\mathbb{H}^{\varepsilon}_{1,u_k}$, $1\leq k\leq s$, separate asymptotically (i.e.,
  there exists a~consistent test procedure for testing $\mathbb{H}_{0,u_k}$ versus $\mathbb{H}^{\varepsilon}_{1,u_k}$, $1\leq k\leq s$)
  when
  the quantity $\min_{1\leq k\leq s} r_{\varepsilon,k}$ is not ``too small''.
The sharp selection boundary in the problem at hand could be described in terms of $\min_{1\leq k\leq s} r_{\varepsilon,k}$, however, the construction of an adaptive selection procedure in Section~\ref{CAS} suggests that it is more natural to describe the sharp selection boundary in terms of $\min_{1\leq k\leq s} a_{\varepsilon,u_k}(r_{\varepsilon,k})/\sqrt{\log {d\choose k}}$, where $a_{\varepsilon,u_k}(r_{\varepsilon,k})$ is defined below by (\ref{def:a2}).

Let $r_\varepsilon=\{r_{\varepsilon,k}, 1\leq k\leq s\}$ be the same family  of collections as in (\ref{hypotheses0}).  Define the class of sparse multivariate functions of our interest by
\begin{multline*}
	{\cal F}^{\beta,\sigma}_{s,d}(r_{\varepsilon})=\Bigg\{f: f(\tb)=\sum_{k=1}^s\sum_{u_k\in {\cal U}_{k,d}} \eta_{u_k} f_{u_k}(\tb_{u_k}), f_{u_k}\;\mbox{satisfies}\;(\ref{orthcon2}),\\
	 f_{u_k}\in \mathring{\cal F}_{\cb_{u_k}}(r_{\varepsilon,k}), u_k \in {\cal U}_{k,d}, \etab_k=(\eta_{u_k})_{u_k\in {\cal U}_{k,d}}\in {H}^k_{\beta,d}, 1\leq k\leq s \Bigg\},
\end{multline*}
where the dependence of ${\cal F}^{\beta,\sigma}_{s,d}(r_{\varepsilon})$ on the smoothness parameter $\sigma$ is hidden in the coefficients $\cb_{u_k}=(c_{\lb})_{\lb\in\mathring{\mathbb{Z}}_{u_k}}$ defined in (\ref{def:c_lb}).
  First, we find the sharp selection boundary that allows us to verify whether the components of a signal $f\in {\cal F}^{\beta,\sigma}_{s,d}(r_{\varepsilon})$
 are all selectable.
  Next, under the assumption that the components $f_{u_k}$, $u_k\in {\cal U}_{k,d},$ $1\leq k\leq s$, are selectable,
 we construct a selector $\boldsymbol{\hat{\eta}}=(\boldsymbol{\hat{\eta}}_k=(\hat{\eta}_{u_k},u_k\in {\cal U}_{k,d}),1\leq k\leq s)$ that recovers exactly the sparsity pattern of $f$ in model
  (\ref{model1})--(\ref{SC1}) and (\ref{Fu}) in the sense that for all $\beta\in(0,1)$ and $\sigma>0$
 \begin{equation}\label{r1}
			\limsup_{\varepsilon\to 0}	\sup_{f\in {\cal F}_{s,d}^{\beta,\sigma}({r_{\varepsilon}})} \E{\varepsilon,f} |\boldsymbol{\hat{\eta}}-\boldsymbol{\eta} | = 0,
		\end{equation}
where ${\operatorname{E}}_{\varepsilon,f}$ is the expectation with respect to the measure ${\operatorname{P}}_{\varepsilon,f} $, and the expression
 $${\operatorname{E}}_{\varepsilon,f} | \boldsymbol{\hat{\eta}}-\boldsymbol{\eta}| ={\operatorname{E}}_{\varepsilon,f}\left(\sum_{k=1}^s\sum_{u_k\in{\cal U}_{k,d}}|\hat{\eta}_{u_k}-{\eta}_{u_k}|\right)$$ is the \textit{Hamming risk} of $\boldsymbol{\hat{\eta}}$ as an estimator of $\etab$.
 The Hamming risk represents the average Hamming loss $|\boldsymbol{\hat{\eta}}-\boldsymbol{\eta}|:=\sum_{k=1}^s\sum_{u_k\in{\cal U}_{k,d}}|\hat{\eta}_{u_k}-\eta_{u_k}|$,
 which counts the number of positions at which $\boldsymbol{\hat{\eta}}$ and $\boldsymbol{{\eta}}$ differ.
  Finally, we show that if at least one of the components $f_{u_k}$ is not selectable,
 then exact recovery of the sparsity pattern of $f\in {\cal F}^{\beta,\sigma}_{s,d}(r_{\varepsilon})$ is impossible.

A related measure of risk that is often used in the literature on variable selection is the
\textit{probability of wrong recovery}, $\operatorname{P}_{\varepsilon,f}(\boldsymbol{\hat{\eta}}\neq \etab)$, that determines
the chance that a selector $\boldsymbol{\hat{\eta}}$ will not agree
with $\etab$. Most literature on variable selection in high dimensions  focuses
on constructing selectors such that the probability of wrong recovery is close to zero
in some asymptotic sense (see, for example, Wainwright (2009), Wasserman and Roeder (2009),
Zhang (2010), Zhao and Yu (2006)).
Although the probability of wrong recovery has been used extensively by other authors, the Hamming risk
is a more general measure of risk since, by means of Markov's inequality,
\begin{gather*}
\operatorname{P}_{\varepsilon,f}(\boldsymbol{\hat{\eta}}\neq \etab)=\operatorname{P}_{\varepsilon,f}(|\boldsymbol{\hat{\eta}}-\boldsymbol{\eta}|\geq 1)
\leq {\operatorname{E}}_{\varepsilon,f} | \boldsymbol{\hat{\eta}}-\boldsymbol{\eta} |.
\end{gather*}
In view of the last inequality, if the Hamming risk tends to zero, the probability of wrong recovery is guaranteed to tend to zero.
For this reason, we only consider the Hamming risk as a measure of selector performance.

\subsection{Sparse Gaussian sequence space model}\label{SGSSM}
We shall study the  recovery problem at hand in the sequence space of Fourier coefficients of $f$
and name it the variable selection problem. A Gaussian sequence space model is equivalent to the corresponding
Gaussian white noise model but is more convenient to deal with since it is written in terms of the Fourier coefficients.

Let $\{\phi_{\lb}(\tb)\}_{\lb\in \mathbb{Z}^d}$ be the orthonormal basis of $L_2^d$ as in (\ref{OB}), and
let $\theta_{\lb}(u_k)=(f_{u_k},\phi_{\lb})_{L_2^d}$ be the $\lb$th Fourier coefficient of $f_{u_k}$ for
$\lb \in 	 \mathring{\mathbb{Z}}_{u_k}$, $u_k\in{\cal U}_{k,d}$, $1\leq k\leq d$.
Denote by  $X_\lb=\mathbb{X}_\varepsilon(\phi_{\lb})=\int_{[0,1]^d}\phi_{\lb}(\tb)\, dX_{\varepsilon}(\tb)$,
$\lb \in  \mathring{\mathbb{Z}}_{u_k}$, $u_k \in \mathcal{U}_{k,d}$, $1 \leq k \leq s$, the $\lb$th empirical Fourier coefficient.
Then, in view of (\ref{orthcon2}) and (\ref{decomp0}), the sequence space model that corresponds to model  (\ref{model1})--(\ref{SC1}) and (\ref{Fu}) takes the form
\begin{eqnarray}\label{model2_s}
X_{\lb}&=&(f,\phi_{\lb})_2+ \varepsilon\mathbb{W}(\phi_{\lb})  =\sum_{k=1}^s\sum_{v_k\in {\cal U}_{k,d}}\eta_{v_k}(f_{v_k},\phi_\lb)_{2}+\varepsilon\mathbb{W}(\phi_{\lb})\nonumber\\ &=&
\eta_{u_k}\theta_{\lb}(u_k)+\varepsilon\xi_\lb, \quad \lb \in \mathring{\mathbb{Z}}_{u_k},\quad  u_k \in \mathcal{U}_{k,d}, \quad 1 \leq k \leq s,
\end{eqnarray}
where $\etab_k=(\eta_{u_k}, u_k\in {\cal U}_{k,d})\in {H}^{k}_{\beta,d}$,
 $ \xi_\lb=\mathbb{W}(\phi_{\lb})=\int_{[0,1]^d}\phi_{\lb}(\tb)\, dW(\tb)$
are iid standard normal random variables for all $\lb \in \mathring{\mathbb{Z}}_{u_k},$  $u_k\in {\cal U}_{k,d}$, $ 1\leq k\leq s,$
 $\thetab_{u_k}=(\theta_{\lb}(u_k),{\lb\in \mathring{\mathbb{Z}}_{u_k}})$ consists of the
  Fourier coefficients
 $\theta_\lb(u_k)=(f_{u_k},\phi_\lb)_{L_2^d}$ and belongs to the ellipsoid
$$	{\Theta}_{\cb_{u_k}}=\bigg \{\thetab_{u_k}=(\theta_{\lb}(u_k),{\lb\in \mathring{\mathbb{Z}}_{u_k}})\in \l2(\mathbb{Z}^d):
	\sum_{\lb \in \mathring{\mathbb{Z}}_{u_k}}\theta_{\lb}^2(u_k)c_{\lb}^2\leq 1 \bigg \}.$$
Next, for $u_k \in {\cal U}_{k,d}$ and $r_{\varepsilon,k}>0$, define
\begin{gather}\label{Theta}	
\mathring{\Theta}_{\cb_{u_k}}(r_{\varepsilon,k})=\bigg \{\thetab_{u_k}\in {\Theta}_{\cb_{u_k}}:
	 \sum_{\lb \in \mathring{\mathbb{Z}}_{u_k}}\theta_{\lb}^2(u_k)\geq r_{\varepsilon,k}^2 \bigg \},
\end{gather}
and note that $\mathring{\Theta}_{\cb_{u_k}}(r_{\varepsilon,k})=\emptyset$ when $r_{\varepsilon,k}>1/c_{\varepsilon,0},$ where,
recalling (\ref{def:c_lb}), $c_{\varepsilon,0}:=\inf_{\lb\in \mathring{\mathbb{Z}}_u}c_{\lb}=(2\pi)^{\sigma}k^{\sigma/2}$.
Therefore, we are only interested in the case when $r_{\varepsilon,k}\in (0,(2\pi)^{-\sigma}k^{-\sigma/2}).$

In the sequence space of Fourier coefficients, the hypothesis testing problem  (\ref{HT1}) becomes that of testing
  \begin{gather}\label{HT2}
 	{H}_{0,u_k} :  \thetab_{u_k}=0 \quad\mbox{vs.}\quad {H}^{\varepsilon}_{1,u_k}:  \thetab_{u_k}\in \mathring{\Theta}_{\cb_{u_k}}(r_{\varepsilon,k}),\quad
 \mbox{for some}\;\; u_k\in{\cal U}_{k,d},\;\;
   1\leq k\leq s.
 \end{gather}
The collection of hypothesis testing problems in (\ref{HT2}) will be used to establish conditions for the possibility of exact variable selection in
model (\ref{model2_s}), and to design an adaptive selection procedure that achieves this type of selection.

For the family of collections $r_\varepsilon=\{r_{\varepsilon,k},\,1\leq k\leq s\}$, $r_{\varepsilon,k}>0$, define the set
\begin{gather*}
\Theta_{s,d}^{\sigma} (r_{\varepsilon})=\{\thetab=(\thetab_1,\ldots,\thetab_s): \thetab_k\in
\mathring{\Theta}^{\sigma}_{k,d}(r_{\varepsilon,k}),1\leq k\leq s\},
\end{gather*}
where
\begin{gather*}
\mathring{\Theta}^{\sigma}_{k,d}(r_{\varepsilon,k})=\{\thetab_k=(\thetab_{u_k},u_k\in{\cal U}_{k,d}):
\thetab_{u_k}\in\mathring{\Theta}_{\cb_{u_k}}(r_{\varepsilon,k}),\, u_k\in {\cal U}_{k,d}\},\quad 1\leq k\leq s.
\end{gather*}
For $u_k\in {\cal U}_{k,d}$, $1\leq k\leq s$, denote ${\boldsymbol X}_{\!u_k}=\{X_{\lb},\lb\in\mathring{\mathbb{Z}}_{u_k}\}$ and let $\boldsymbol{\hat{\eta}}_k=(\hat{\eta}_{u_k},u_k \in {\cal U}_{k,d} )$, where $\hat{\eta}_{u_k} =\hat{\eta}_{u_k}({\boldsymbol X}_{\!u_k})\in\{0,1\}$, be an estimator of $\etab_k=({\eta}_{u_k},u_k \in {\cal U}_{k,d} )\in { H}_{\beta,d}^{k}.$
We have previously decided to call an aggregate estimator $\boldsymbol{\hat{\eta}}=(\boldsymbol{\hat{\eta}}_1,\ldots,\boldsymbol{\hat{\eta}}_s)$
for $\etab=(\etab_1,\ldots,\etab_s)\in  \mathcal{H}^s_{\beta,d}$ a \textit{selector}.
When dealing with the problem of identifying nonzero $\eta_{u_k}$ in model (\ref{model2_s}),
the \textit{maximum Hamming risk} of a selector $\boldsymbol{\hat{\eta}}$ can be expressed as
$$R_{\varepsilon,s}(\boldsymbol{\hat{\eta}}):=\sup\limits_{\boldsymbol{\eta} \in \mathcal{H}^s_{\beta,d}} \sup\limits_{\thetab\in  \Theta_{s,d}^{\sigma}(r_{\varepsilon})} \E{\thetab,\etab} | \boldsymbol{\hat{\eta}}-\boldsymbol{\eta}|,$$
where  ${\operatorname{E}}_{\thetab,\etab}$ is the expectation with respect
 to the joint distribution of $\boldsymbol X_{\!u_k}=\{X_{\lb},\lb\in\mathring{\mathbb{Z}}_{u_k}\}$, $u_k\in{\cal U}_{k,d},$ $1\leq k\leq s$,
  in model (\ref{model2_s}), and
 $${\operatorname{E}}_{\thetab,\etab} | \boldsymbol{\hat{\eta}}-\boldsymbol{\eta} | ={\rm E}_{\thetab,\etab}\left(\sum_{k=1}^s\sum_{u_k\in{\cal U}_{k,d}}|\hat{\eta}_{u_k}-{\eta}_{u_k}|\right)$$ is the \textit{Hamming risk} of $\boldsymbol{\hat{\eta}}$.
A ``good'' selector is the one that consistently chooses the nonzero $\eta_{u_k}$ as $\varepsilon\to 0$. With this in mind,
an \textit{exact selector} $\hat{\etab}$  is defined such that $\limsup_{\varepsilon\to 0}R_{\varepsilon,s}(\boldsymbol{\hat{\eta}})=0$ for all $\beta\in(0,1)$ and $\sigma>0$.
If such a selector exists, it is said that exact selection in model  (\ref{model2_s}) is possible.

\subsection{Summary of the main results}
Consider  the problem of
exact identification of the nonzero index variables $\eta_{u_k}$
in the sequence space model  (\ref{model2_s}), and recall the family of
collections $r_\varepsilon=\{r_{\varepsilon,k},1\leq k\leq s\}$ that determines
the sets $ \mathring{\Theta}_{\cb_{u_k}}(r_{\varepsilon,k}) $, $u_k\in{\cal U}_{k,d}$, $1\leq k\leq s$,
as defined in (\ref{Theta}).
The main contribution of this work is the derivation of the \textit{sharp selection boundary}, stated in terms of  $r_\varepsilon=\{r_{\varepsilon,k},1\leq k\leq s\}$,
which defines when exact selection of the nonzero index variables $\eta_{u_k}$ in model (\ref{model2_s}) is possible and when this type of selection is impossible.
 This boundary is determined by inequalities {\rm (\ref{cond1})} and  {\rm (\ref{cond:inf2})} bellow.
Then, we construct
an adaptive (independent of the sparsity index $\beta$) estimator $\boldsymbol{\hat{\eta}}$ of $\etab\in {\cal H}^{s}_{\beta,d}$ attaining this boundary with the property
\begin{equation}\label{r11}
			\limsup_{\varepsilon\to 0}	\sup_{\boldsymbol{\eta} \in \mathcal{H}^s_{\beta,d}} \sup_{\thetab\in  \Theta_{s,d}^{\sigma}(r_{\varepsilon})} \E{\thetab,\etab} | \boldsymbol{\hat{\eta}}-\boldsymbol{\eta}| = 0,
		\end{equation}
 which holds for all $\beta\in(0,1)$ and  $\sigma>0$.
  We also show that if the family of collections
$r_\varepsilon=\{r_{\varepsilon,k},1\leq k\leq s\}$ falls below the selection boundary,
one has, for $\beta\in(0,1)$ and $\sigma>0$,
\begin{equation}\label{r22}
		\liminf_{\varepsilon \to 0}\inf_{\tilde{\etab}} \sup_{\boldsymbol{\eta} \in  \mathcal{H}^s_{\beta,d}} \sup_{\thetab\in \Theta_{s,d}^{\sigma}(r_{\varepsilon})} \E{\thetab,\etab} | \boldsymbol{\tilde{\eta}}-\boldsymbol{\eta}|>0,
	\end{equation}
where the infimum is taken over all selectors  $\boldsymbol{\tilde{\eta}}$ of $\boldsymbol{\eta} \in \mathcal{H}^s_{\beta,d}$ in model
(\ref{model2_s}), and so, exact  selection of the nonzero index variables $\eta_{u_k}$ in model (\ref{model2_s}) is impossible.
 This  implies that the proposed selector is optimal (in the  minimax hypothesis testing sense).

 The limiting relations (\ref{r11}) and (\ref{r22}) will be referred to as the \textit{upper bound} on the maximum Hamming risk of $\boldsymbol{\hat{\eta}}$
 and the \textit{lower bound} on the minimax Hamming risk, respectively.

\subsection{Previous results}
In order to better situate our work in current literature, consider the problem of
variable selection in high-dimensional regression.
 This problem, in its simplest form, consists of identifying active (nonzero) coordinates of a nonrandom vector $\boldsymbol{x}=(x_1,\ldots,x_m)^{\top} \in \mathbb{R}^m$ observed with additive noise:
\begin{equation}
	\boldsymbol{Y} = A \boldsymbol{x} + \boldsymbol{Z},
	\label{sparse_recovery_prob}
\end{equation}
where $\boldsymbol{Y}=(Y_1,\ldots,Y_n)^{\top}$ is a response vector,  $\boldsymbol{Z}=(Z_1,\ldots,Z_n)^{\top}$ is a zero mean noise vector
with some (usually diagonal) covariance matrix,
$A$ is an $n\times m$ design matrix that is either deterministic (and often equal to identity) or has some known distribution,
and the dimension $m$ is much larger than the sample size $n$.
To ensure that the recovery of the ``signal'' vector $\boldsymbol{x}$ in model (\ref{sparse_recovery_prob}) is feasible,
the vector $\boldsymbol{x}$ is assumed to be sparse.
Depending on the context, the recovery of the sparsity pattern of $\boldsymbol{x}$ in model  (\ref{sparse_recovery_prob}) is termed in the literature as variables selection, subset selection, support recovery, or multiple testing (see, for example, Wasserman and Roeder (2009), Belitser and Nurushev (2018), Butucea et al. (2018), Gao and Stove (2020), and references wherein).
 However, the goal is always the same -- to correctly detect  the support $\boldsymbol{S} = \{j: x_j \neq 0 \}$
of $\boldsymbol{x}$. In addition to finding an estimator $\boldsymbol{\hat{S}}_n(\boldsymbol{Y})$ of $\boldsymbol{S}$ that recovers the positions of nonzero entries of $\boldsymbol{x}$,
which is often called a selector or decoder of $\boldsymbol{S}$, it is also important to establish the minimal conditions under which the support set can be effectively estimated.
The quality of an estimator of $\boldsymbol{S}$ can be measured by various metrics, including the Hamming loss, defined as the number of mismatches between the estimated and true support set,
the probability of support recovery $\operatorname{P}(\boldsymbol{\hat{S}}_n(\boldsymbol{Y}) = \boldsymbol{S})$, or the misdetection probability, defined as the probability that the number of misdetected components equals zero. The literature on support recovery in the high-dimensional linear model (\ref{sparse_recovery_prob}) under sparsity assumption is very rich and its complete overview would fall far beyond the format of this paper. So, we provide only an outline of some results that are relevant to our work.

Two most common approaches to variable selection are via penalization and thresholding. Penalization methods are suitable for situations where the dimension $m$ is not very large and their core is imposing a proper penalty on the number of selected variables. The penalty can be given simply by $l_1$-norm of $\boldsymbol{x}$ or based on some information criterion (e.g., AIC, BIC)  or even more complex (e.g., concave penalty as in Zhang (2010)). One of the most popular penalization methods, lasso, uses the $l_1$-norm of $\boldsymbol{x}$ as a penalty term (see, for example,  Zhao and Yu (2006), Wainwright (2009), Wasserman and Roeder (2009)). Lasso is efficient for relatively large $m$, but suffers from nonnegligible bias. For very large $m$,  penalization becomes computationally demanding and methods based on univariate thresholding, also called marginal regression, come into play (see
Genovese et al. (2012) and Ndaoud (2020)). An estimator of $\boldsymbol{S}$ provided by a thresholding method has a simple form $\boldsymbol{\hat{S}}_n (\boldsymbol{Y})= \{j: |Y_j| > t(\boldsymbol{Y})\}$ and,
with a suitable choice of  $t(\boldsymbol{Y})$, it is very fast and efficient.
The threshold $t(\boldsymbol{Y})$ can be set up in many different ways  and in case of some iterative methods, it can be even updated at each step (see Ndaoud (2020)).
The selection procedure introduced in this work also belongs to the class of thresholding methods.

The problem of variable selection in sparse nonparametric regression is far more challenging
as compared to its parametric analogue. Therefore, the related results obtained thus
far are not as numerous and diverse as those obtained in parametric settings.
In a simpler setup, when instead of decomposition (\ref{fun_sub_anova22}) in model (\ref{model1})--(\ref{SC1}) and (\ref{Fu})
the signal $f(\tb)$ is assumed to be a sum of functions each depending on $k$ variables ($1\leq k\leq d$) as in (\ref{f}),
the problem of identifying the active components of $f(\tb)$ was studied in Ingster and Stepanova (2014) and Butucea and Stepanova (2017) for $k=1$, and
in Stepanova and Turcicova (2025) for $1\leq k\leq d$. In these articles, conditions under which exact and almost full identification of the nonzero components
of $f(\tb)$ as in (\ref{f})  is possible and impossible were established,
and adaptive selection procedures that, under some model assumptions, identify exactly and almost fully all nonzero $k$-variate components $f_{u_k}$
were proposed. For more details, we refer to the discussion of Section~\ref{Discussion}.

\section{Construction of adaptive selector}\label{CAS}
Return to the problem of testing $H_{0,u_k}$ against $H^{\varepsilon}_{1,u_k}$ for some $u_k\in{\cal U}_{k,d}$, $1\leq k\leq s$, as in (\ref{HT2}) and note the following facts from Ingster and Suslina (2003, 2005) and Stepanova and Turcicova (2025).
For $u_k\in {\cal U}_{k,d}$ and $r_{\varepsilon,k}>0$, $1\leq k\leq s$, consider the value $a^2_{\varepsilon,u_k} (r_{\varepsilon,k})$ of the following extreme problem:
	\begin{equation}
		a^2_{\varepsilon,u_k} (r_{\varepsilon,k})= \frac{1}{2 \varepsilon^4} \inf_{\thetab_{u_k} \in \mathring{\Theta}_{\cb_{u_k}}(r_{\varepsilon,k})} \sum_{\lb\in \mathring{\mathbb{Z}}_{u_k}} \theta^4_\lb (u_k),
		\label{def:a2}
	\end{equation}
where, as before,
$$	\mathring{\Theta}_{\cb_{u_k}}(r_{\varepsilon,k})=\bigg \{\thetab_{u_k}=(\theta_{\lb}(u_k),{\lb\in \mathring{\mathbb{Z}}_{u_k}})\in \l2(\mathbb{Z}^d):
	\sum_{\lb \in \mathring{\mathbb{Z}}_{u_k}}\theta_{\lb}^2(u_k)c_{\lb}^2\leq 1,\;	 \sum_{\lb \in \mathring{\mathbb{Z}}_{u_k}}\theta_{\lb}^2(u_k)\geq r_{\varepsilon,k}^2 \bigg \}.$$
 The function $a_{\varepsilon,u_k} (r_{\varepsilon,k})$ plays a key role in the minimax theory of hypothesis testing.
It controls the minimax total error probability and is used
to set up a cut-off point of the asymptotically minimax test procedure
in the problem of testing ${H}_{0,u_k}:\thetab_{u_k}=\boldsymbol{0}$ versus ${H}^{\varepsilon}_{1,u_k}:\thetab_{u_k} \in \mathring{\Theta}_{\cb_{u_k}}(r_{\varepsilon,k})$
(see, for example, Theorem 2 of Ingster and Suslina (2005)).

Note that $a_{\varepsilon,u_k} (r_{\varepsilon,k})$ is a non-decreasing function of its argument that possesses a kind of ``continuity'' property
 in the sense that if $r_{\varepsilon,k}/r^{\prime}_{\varepsilon,k}$ is close to 1, then $a_{\varepsilon,u_k} (r_{\varepsilon,k})/a_{\varepsilon,u_k} (r^{\prime}_{\varepsilon,k})$ is also close to 1 and vise versa.
 Namely, for all $\gamma >0$, there exists $\varepsilon^*>0$ and $\delta^* >0$ such that (see Section 5.2.3 of Ingster and Suslina (2003))
\begin{equation*}
	a_{\varepsilon,u_k} (r_{\varepsilon,k}) \leq a_{\varepsilon,u_k} ((1+\delta)r_{\varepsilon,k}) \leq (1+\gamma) a_{\varepsilon,u_k} (r_{\varepsilon,k}), \quad \forall\, \varepsilon \in (0,\varepsilon^*), \forall\, \delta \in (0,\delta^*).
\end{equation*}

Suppressing for brevity the dependence on $u_k$, denote by $\{\theta^*_\lb(r_{\varepsilon,k}),\lb\in \mathring{\mathbb{Z}}_{u_k}\}$ the extreme sequence in (\ref{def:a2}), that is,
\begin{gather}
a^2_{\varepsilon,u_k} (r_{\varepsilon,k})=   \frac{1}{2 \varepsilon^4} \sum_{\lb\in \mathring{\mathbb{Z}}_{u_k}} \left(\theta^*_\lb (r_{\varepsilon,k}) \right)^4.
\label{def:a2_theta_star}
\end{gather}
The sharp asymptotics for $\theta^*_\lb (r_{\varepsilon,k})$ when $k$ is fixed and when $k=k_\varepsilon\to \infty$, $k=o(d)$ as $\varepsilon\to 0$
are known (see Section 2 of Stepanova  and Turcicova (2025)). For fixed $k$, as $\varepsilon\to 0$,
 \begin{gather}\label{thetal}
	(\theta^*_\lb (r_{\varepsilon,k}))^2 \sim
	\frac{r_{\varepsilon,k}^{2+{k}/{\sigma}} 2^k \pi^{{k}/{2}}(k+2\sigma)\Gamma\left(1+{k}/{2}\right) }{2\sigma\left(1+4\sigma/k  \right)^{{k}/{(2\sigma)}}}
\left(1-\left(\sum\nolimits_{j=1}^d (2\pi l_j)^2  \right)^{\sigma}\frac{r^2_{\varepsilon,k} }{(1+4\sigma/k)}  \right)_+,
\end{gather}	
where $x_+=\max(x,0)$.
For $k=k_\varepsilon\to \infty$, $k=o(d)$ as $\varepsilon\to 0$, the above relation
 can be modified by using $\Gamma(x+1)\sim \sqrt{2\pi}x^{x+1/2}e^{-x}$ and $(1+1/x)^x\sim e$ as $x\to \infty$.
Expression (\ref{thetal}) implies that
  \begin{gather*}
\#\{\lb \in \mathring{\mathbb{Z}}_{{u_k}}: \theta^*_\lb (r_{\varepsilon,k})\neq 0\}
=\#\left\{\lb \in \mathring{\mathbb{Z}}_{{u_k}}:\left(\sum_{j=1}^d l_j^2  \right)^{1/2}< (1+4\sigma/k)^{1/(2\sigma)}/(2\pi r_{\varepsilon,k}^{1/\sigma})\right\}
\end{gather*}
is of the same order of magnitude as the integer-volume of the $\l2$-ball of radius $(1+4\sigma/k)^{1/(2\sigma)}/(2\pi r_{\varepsilon,k}^{1/\sigma})$ in $\mathbb{R}^k$,
which is  $O(r_{\varepsilon,k}^{-{k}/{\sigma}})$,
and hence
 \begin{gather}\label{NL}
 \#\{\lb\in \mathring{\mathbb{Z}}_{u_k}: \theta^*_\lb (r_{\varepsilon,k})\neq 0\}  ={O}( r_{\varepsilon,k}^{-k/\sigma}),\quad \varepsilon\to 0,
 \end{gather}
 when $k$ is fixed ($1\leq k\leq d$) and when  $k=k_\varepsilon\to \infty$, $k=o(d)$ as $\varepsilon\to 0$.

	It is known  that for a fixed $k$, the sharp asymptotics of $a_{\varepsilon,k}(r_{\varepsilon,k})$ as $\varepsilon\to 0$  are given by
(see, for example, Theorem 4 of Ingster and Suslina (2005))
	\begin{equation}\label{aek}
		a_{\varepsilon,u_k}(r_{\varepsilon,k})\sim C(\sigma,k)r_{\varepsilon,k}^{2+{k}/{(2\sigma)}}\varepsilon^{-2},\quad
		C^2(\sigma,k)=\frac{\pi^k (1+{2\sigma}/{k})\Gamma(1+k/2)}{(1+4\sigma/k)^{1+{k}/{(2\sigma)}}\Gamma^k(3/2)},
	\end{equation}
	and when $k=k_\varepsilon\to \infty$, $k=o(d)$ as $\varepsilon\to 0$, they are given by (see, for example, relation (39) in Ingster and Suslina (2005))
	\begin{equation}\label{aek1}
		a_{\varepsilon,u_k}(r_{\varepsilon,k})\sim \left({2\pi k}/{e}\right)^{{k}/{4}}e^{-1}(\pi k)^{1/4}r_{\varepsilon,k}^{2+{k}/{(2\sigma)}}\varepsilon^{-2}.
	\end{equation}
The \textit{sharp selection boundary} in the problem under study, which is given
in Section~\ref{sec:main} by Theorems 1 and 3 for a fixed $s$, and by Theorems 2 and 4 for $s=s_\varepsilon\to \infty$, $s=o(d)$ as $\varepsilon\to \infty$,
is described in terms of $\min_{1\leq k\leq s}a_{\varepsilon,u_k}(r_{\varepsilon,k})/\sqrt{\log {d\choose k}}$.

If the sparsity index $\beta$ was known, we could construct an exact selector of $\boldsymbol{\eta}=(\etab_1,\ldots,\etab_s) \in \mathcal{H}^s_{\beta,d}$ as follows.
First, determine a positive family of collections $r^*_\varepsilon=\{r^*_{\varepsilon,k},1\leq k\leq s\}$ by, cf. condition~(\ref{cond1}) in Theorem \ref{Theorem1} below,
\begin{gather} \label{def:r_star_ekm}
	a_{\varepsilon,u_k} (r^*_{\varepsilon,k}) = (1+\sqrt{1-\beta})\sqrt{2\log {d\choose k}},\quad \mbox{for some}\; u_k\in{\cal U}_{k,d},\; 1\leq k\leq s,
\end{gather}
and consider weighted $\chi^2$-type statistics
\begin{equation*}
	S_{u_k}(\beta) = \sum_{\lb \in \mathring{\mathbb{Z}}_{u_k}} \omega_\lb (r^*_{\varepsilon,k}) \left( \left( {X_\lb}/{\varepsilon}  \right)^2 -1 \right),
\quad u_k \in {\cal U}_{k,d},\quad 1\leq k\leq s,
	\label{def:S_u}
\end{equation*}
where  the weight function $\omega_\lb (r_{\varepsilon,k}) $, $r_{\varepsilon,k}>0$, is given by
\begin{equation*}
	\omega_\lb (r_{\varepsilon,k}) = \frac{1}{2 \varepsilon^2} \frac{(\theta^*_\lb (r_{\varepsilon,k}))^2}{a_{\varepsilon,u_k}(r_{\varepsilon,k})}, \quad \lb \in \mathring{\mathbb{Z}}_{u_k},\quad u_k\in{\cal U}_{k,d},\;1\leq k\leq s.
	\label{def:omega}
\end{equation*}
 It follows from (\ref{def:a2}) that
$	\sum_{\lb \in \mathring{\mathbb{Z}}_{u_k}} \omega_\lb^2 (r_{\varepsilon,k}) ={1}/{2}$ {for all} $r_{\varepsilon,k} > 0.$
	The statistic $S_{u_k}(\beta)$ depends on $\beta$ through the weights $\omega_\lb (r^*_{\varepsilon,k})$, $\lb \in \mathring{\mathbb{Z}}_{u_k}$,
and, in view of (\ref{NL}), consists of $O((r^*_{\varepsilon,k})^{-k/\sigma})$ nonzero terms.
Next, for $u_k\in {\cal U}_{k,d}$, $1\leq k\leq s$, consider the data ${\boldsymbol X}_{\!u_k}=\{X_{\lb},\lb\in\mathring{\mathbb{Z}}_{u_k}\}$
observed in model (\ref{model2_s}) and let $\boldsymbol{\check{\eta}}_k(\beta)=(\check{\eta}_{u_k}(\beta),u_k \in {\cal U}_{k,d} )$, where $\check{\eta}_{u_k}(\beta)$
depends on ${\boldsymbol X}_{\!u_k}$  and takes on its values in $\{0,1\}$, be an estimator of $\etab_k=({\eta}_{u_k},u_k \in {\cal U}_{k,d} )\in { H}_{\beta,d}^{k}$ defined by
\begin{gather}\label{eta_check}
\check{\eta}_{u_k}(\beta) =\ind{ S_{u_k}(\beta)>\sqrt{(2+\epsilon)\log{d\choose k}} },\quad u_k\in {\cal U}_{k,d},\quad 1\leq k\leq s,
\end{gather}
where $\epsilon>0$ depends on $\varepsilon$ and tends to zero according to either (\ref{delta}) or (\ref{delta1}),
depending on whether $s$ is fixed or $s\to \infty$, $s=o(d)$ as $\varepsilon\to 0.$
Then, under the assumptions of  Theorem \ref{Theorem1} in Section \ref{UpBnd} if $s$ is fixed, and under the assumptions of Theorem \ref{Theorem2}
in Section \ref{UpBnd} if
$s\to \infty$, $s=o(d)$ as $\varepsilon\to 0 $,
 an aggregate selector $\boldsymbol{\check{\etab}}(\beta)=(\boldsymbol{\check{\etab}}_k(\beta)=(\check{\eta}_{u_k}(\beta),u_k\in {\cal U}_{k,d}),1\leq k\leq s)$, with  $\check{\eta}_{u_k}(\beta)$ as in (\ref{eta_check}), is an exact selector of $\boldsymbol{{\eta}}\in {\cal H}^s_{\beta,d}$. The proof of this claim is omitted because
it largely goes along the lines of the proofs of Theorems \ref{Theorem1} and \ref{Theorem2} in Section \ref{UpBnd}, only shorter and easier.

However, if $\beta$ is unknown, the selector $\boldsymbol{\check{\etab}}(\beta)$ is not applicable.
The construction of an adaptive exact selector working for all values of $\beta\in (0,1)$ and attaining the sharp selection boundary utilizes the Bonferonni idea and is
similar to that in Section 2 of Stepanova and Turcicova (2025).
	For a~given~$k$ ($1\leq k\leq s$), consider a grid of equidistant points on the interval (0,1):
	\begin{equation*}
		\beta_{m,k} = m \Delta_k,\quad m=1,\ldots,M_k,
	\end{equation*}
	where $\Delta_k>0$ is a small number {that depends on $\varepsilon$, $M_k = [ 1/\Delta_k] \in \mathbb{N}$, and
		$\log M_k = o \left( \log {d\choose k} \right)$, $1\leq k\leq s$, as
		$\varepsilon \to 0$.}
	Next, for $m=1,\ldots,M_k$, $1\leq k\leq s$, let $r^*_{\varepsilon,k,m} >0$ be the solution of the equation
	\begin{equation*}
		a_{\varepsilon,u_k} (r^*_{\varepsilon,k,m}) = (1+\sqrt{1-\beta_{m,k}}) \sqrt{2 \log {d\choose k} }.
	\end{equation*}
It can be shown that for $m=1,\ldots,M_k$, $1\leq k\leq s$,
 $r^*_{\varepsilon,k,m}=O\left( \left(\varepsilon^4 \log {d\choose k}   \right)^{\sigma/(4\sigma+k)}\right)$ when $k$ is fixed, and
$r^*_{\varepsilon,k,m}=O\left( \left(\varepsilon^4 \log {d\choose k}   \right)^{\sigma/(4\sigma+k)}k^{-\sigma/2}\right)$ when $k=k_\varepsilon\to \infty$,
as $\varepsilon\to 0$  (see relations (45) and (92) in Stepanova and Turcicova (2025)).

For every $u_k \in {\cal U}_{k,d}$, $1\leq k\leq s$, consider the weighted $\chi^2$-type statistics
	\begin{equation*}
		S_{u_k,m} = \sum_{\lb \in \mathring{\mathbb{Z}}_{u_k}} \omega_\lb (r^*_{\varepsilon,k,m}) \left(  ({X_\lb}/{\varepsilon})^2 -1 \right), \quad m=1,\ldots,M_k,
	\end{equation*}
	where
	\begin{equation*}
		\omega_\lb (r_{\varepsilon,k,m}) = \frac{1}{2 \varepsilon^2} \frac{(\theta^*_\lb (r_{\varepsilon,k,m}))^2}{a_{\varepsilon,u_k}(r_{\varepsilon,k,m})}, \quad \lb \in \mathring{\mathbb{Z}}_{u_k}.
	\end{equation*}
Recall that $  \sum_{\lb \in \mathring{\mathbb{Z}}_{u_k}}\omega^2_\lb (r_{\varepsilon,k,m}) =1/2 $ for all $r_{\varepsilon,k,m}>0$, and, in view of (\ref{NL}),
 every statistic $S_{u_k,m}$ consists of $O((r^*_{\varepsilon,k,m})^{-k/\sigma})$ nonzero terms for all $u_k\in {\cal U}_{k,d}$, $1\leq m\leq M_k$, $1\leq k\leq s$.
 It is also known that for all $u_k\in {\cal U}_{k,d}$, $1\leq m\leq M_k$, $1\leq k\leq s$ (see equation (47) of Stepanova and Turcicova (2025)),
 \begin{equation*}
	\max_{\lb \in \mathring{\mathbb{Z}}_{u_k}} \omega_\lb(r_{\varepsilon,k,m}^*) \asymp \left( \varepsilon \left\{ \log {d\choose k} \right\}^{{1}/{4}} \right)^{{2k}/{(4 \sigma +k)}}.
\end{equation*}

The statistics $S_{u_k,m}$, $m=1,\ldots,M_k$, are similar to those that determine the asymptotically minimax tests for testing
${H}_{0,u_k}:\thetab_{u_k}=\boldsymbol{0}$ versus ${H}^{\varepsilon}_{1,u_k}:\thetab_{u_k} \in \mathring{\Theta}_{\cb_{u_k}}(r_{\varepsilon,k})$, $u_k\in {\cal U}_{k,d}$
(see Theorem 2 of Ingster and Suslina (2005)). These statistics are also useful for the purpose of constructing an exact selector,
as demonstrated below.

Before proceeding, let us note some useful properties of the statistics $S_{u_k,m}$, $m=1,\ldots,M_k$, $1\leq k\leq s$.
These properties are crucial for the adaptive selector $\boldsymbol{\hat{\eta}}=(\boldsymbol{\hat{\eta}}_k=(\hat{\eta}_{u_k},u_k \in {\cal U}_{k,d}),1\leq k\leq s)$ introduced below to be exact.
Denote by $\chi^2_{\nu}(\lambda)$ a noncentral chi-square random variable with $\nu$ degrees of freedom and
noncentrality parameter $\lambda$.
The random variables
$(X_\lb/\varepsilon)^2$, $\lb\in \mathring{\mathbb{Z}}_{u_k}$, $1\leq k\leq s$, are independent and $(X_\lb/\varepsilon)^2\sim \chi_1^2((\theta_\lb/\varepsilon)^2)$
where
$\thetab_{u_k}=(\theta_\lb)_{\lb\in \mathring{\mathbb{Z}}_{u_k}} \in \mathring{\Theta}_{\cb_{u_k}}(r_{\varepsilon,k})$, $1\leq k\leq s$.
Therefore, noting that $\operatorname{E}(\chi_\nu^2(\lambda))=\nu+\lambda$, ${\rm var}(\chi_\nu^2(\lambda))=2(\nu+2\lambda)$,
and  taking into account the normalization condition for the sum of squared weights, it holds, under ${H}^{\varepsilon}_{1,u_k}$,  for $m=1,\ldots,M_k$, $1\leq k\leq s$, that
\begin{align*}
	\E{\thetab_{u_k}}(S_{u_k,m})&=\sum_{\lb\in\mathring{\mathbb{Z}}_{u_k}}\omega_\lb (r^*_{\varepsilon,k,m}) (\theta_\lb/\varepsilon)^2, \nonumber \\
	{\rm var}_{\thetab_{u_k}} (S_{u_k,m})& =\sum_{\lb \in \mathring{\mathbb{Z}}_{u_k}} {\rm var}_{\thetab_{u_k}} \left(\omega_\lb (r_{\varepsilon,k,m}^*) (X_\lb/\varepsilon)^2\right)=
	\sum_{\lb \in \mathring{\mathbb{Z}}_{u_k}} \omega^2_\lb (r_{\varepsilon,k,m}^*)\left(2+4(\theta_\lb/\varepsilon)^2\right)\nonumber\\
	&=
	1+4\sum_{\lb \in \mathring{\mathbb{Z}}_{u_k}} \omega^2_\lb (r_{\varepsilon,k,m}^*)(\theta_\lb/\varepsilon)^2
	=  1+{O}\left(\max_{\lb\in  \mathring{\mathbb{Z}}_{u_k}}\omega_\lb(r_{\varepsilon,k,m}^*)\E{\thetab_{u_k}}(S_{u_k,{m}})\right),\quad \varepsilon\to 0.\label{VarS}
\end{align*}
It is clear that, under ${H}_{0,u_k}$, we have, for $m=1,\ldots,M_k$, $1\leq k\leq s$,
\begin{gather*}
	\operatorname{E}_{\boldsymbol{0}}(S_{u_k,m})=0,\quad {\rm{var}}_{\boldsymbol{0}}(S_{u_k,m})=1.
\end{gather*}
Also, let $T =T_{\varepsilon,k}\rightarrow \infty$ as $\varepsilon\to 0$ be such that
\begin{gather}
	T \max_{\lb \in \mathring{\mathbb{Z}}_{u_k}} \omega_\lb (r^*_{\varepsilon,k,m})=o(1), \quad \varepsilon\to 0.\label{lem:ass_on_T}
\end{gather}
Then, for $m=1,\ldots,M_k,$ $1\leq k\leq s$, as $\varepsilon\to 0$
\begin{gather}\label{bound_1}
	\operatorname{P}_{\boldsymbol{0}} \left( S_{u_k,m} > T \right) \leq \exp \left( - \dfrac{T^2}{2}(1+o(1)) \right).
\end{gather}
If, in addition to (\ref{lem:ass_on_T}), one has for $\thetab_{u_k}\in \mathring{\Theta}_{\cb_{u_k}}(r_{\varepsilon,k})$  and $m=1,\ldots,M_k,$ $1\leq k\leq s$,
\begin{gather*}
	\E{\thetab_{u_k}} (S_{u_k,m}) \max_{\lb \in \mathring{\mathbb{Z}}_{u_k}} \omega_\lb (r^*_{\varepsilon,k,m}) =o(1),
	\label{lem:ass_on_ES}
\end{gather*}
then for this $\thetab_{u_k}$ and $m=1,\ldots,M_k,$ $1\leq k\leq s$, as $\varepsilon\to 0$
\begin{gather}\label{bound_2}
	\operatorname{P}_{\thetab_{u_k}} \left( S_{k,m} - \operatorname{E}_{\thetab_{u_k}}
	(S_{k,m}) \leq -T \right) \leq \exp \left( -\frac{T^2}{2}(1+\smallO(1)) \right).
\end{gather}
The proofs of (\ref{bound_1}) and (\ref{bound_2}) are similar to that of Proposition 7.1 in Gayraud and Ingster (2012) and therefore are omitted.

	Now, we propose a selector $\boldsymbol{\hat{\eta}}=(\boldsymbol{\hat{\eta}}_k=(\hat{\eta}_{u_k},u_k \in {\cal U}_{k,d}),1\leq k\leq s)$
 with components,
	cf. relation~(34) in Stepanova and Turcicova (2025),
	\begin{equation}
		\hat{\eta}_{u_k} = \max_{1 \leq m \leq M_k} \hat{\eta}_{u_k,m}, \quad \quad  \hat{\eta}_{u_k,m} = \ind{ S_{u_k,m} > \sqrt{(2+\epsilon) \left( \log {d\choose k} + \log M_k \right)} },
		\label{selector}
	\end{equation}
	where $\epsilon>0$ depends on $\varepsilon$ and satisfies
	\begin{gather}\label{delta}
		\epsilon\to 0\quad\mbox{and}\quad \epsilon\log {d}\to \infty,\quad \mbox{as}\;\; \varepsilon\to 0,
	\end{gather}
	when $s$ is fixed, and
	\begin{gather}\label{delta1}
		\epsilon\to 0,\quad \epsilon\log d\to \infty,\quad \log s=o(\epsilon \log d),\quad \mbox{as}\;\; \varepsilon\to 0,
	\end{gather}
	when $s=s_\varepsilon\to\infty$, $s=o(d)$ as $\varepsilon\to 0$. Note that the requirement $\epsilon\log d\to \infty$ as $\varepsilon\to 0$
implies $\epsilon\log {d\choose k}\to \infty$ as $\varepsilon\to 0$ for
all $1\leq k\leq s$, where $s$ is either fixed or $s=s_\varepsilon\to \infty$, $s=o(d)$.
Likewise, if for $s=s_\varepsilon\to \infty$ such that $s=o(d)$ one has $\log s=o(\epsilon \log d)$ as $\varepsilon\to 0$, then
$\log s=o(\epsilon \log {d\choose k})$ as $\varepsilon\to 0$  for all $1\leq k\leq s$.
 It should be understood that the proposed selector $\boldsymbol{\hat{\eta}}$
depends on $\sigma>0$, and therefore we have a whole class of adaptive selectors indexed by $\sigma$.

The idea behind the aggregate selector $\boldsymbol{\hat{\eta}}=(\boldsymbol{\hat{\eta}}_k=(\hat{\eta}_{u_k},u_k \in {\cal U}_{k,d}),1\leq k\leq s)$
defined for fixed~$s$ by~(\ref{selector}) and (\ref{delta}) and for $s=s_\varepsilon\to \infty$ by (\ref{selector}) and (\ref{delta1}) is a natural one.
For $u_k\in {\cal U}_{k,d}$, $1\leq k\leq s$, $\hat{\eta}_{u_k}$~identifies the component $f_{u_k}$  as active, when at least one
of the statistics $ S_{u_k,m}$, $m=1,\ldots,M_k$, detects it. Therefore,  the probability of falsely non-identifying
$f_{u_k}$ as active by means of $\hat{\eta}_{u_k}$ does not exceed the probability $\operatorname{P}_{\eta_{u_k}}(\hat{\eta}_{u_k,m_{0,k}} = 0)$
with $\eta_{u_k}=1$
for some $\beta_{m_{0,k},k}$ close to the true (but unknown) value of $\beta$, and $\sum_{k=1}^s \sum_{u_k:\eta_{u_k}=1}\operatorname{P}_{\eta_{u_k}}(\hat{\eta}_{u_k,m_{0,k}} = 0)$ is small
 by the choice of the threshold in the definition of $\hat{\eta}_{u_k,m}$.
Also, the probability of falsely identifying $f_{u_k}$ as active is bounded from above by the sum
$\sum_{m=1}^{M_k} \operatorname{P}_{\eta_{u_k}}(\hat{\eta}_{u_k,m} =1)$ with $\eta_{u_k}=0$, and $\sum_{k=1}^s\sum_{u_k:\eta_{u_k}=0}\sum_{m=1}^{M_k} \operatorname{P}_{\eta_{u_k}}(\hat{\eta}_{u_k,m} =1)$
is small by the choice of the threshold  in the definition of $\hat{\eta}_{u_k,m}$. For more details, see the proof of  Theorem \ref{Theorem1}
in Section \ref{UpBnd}.

\section{Main results} \label{sec:main}
We first state the conditions when exact variable selection in model (\ref{model2_s}) is possible and show that the proposed selector
$\boldsymbol{\hat{\eta}}=(\boldsymbol{\hat{\eta}}_k=(\hat{\eta}_{u_k},u_k \in {\cal U}_{k,d}),1\leq k\leq s)$
  achieves this type of selection. Then, we demonstrate
that our selector is the best possible in the asymptotically minimax sense.
In the statements of Theorems 1 to 4 below, $u_k$ is an arbitrary element of ${\cal U}_{k,d}$ for $1\leq k\leq s$.

\subsection{Upper bound on the maximum risk of $\boldsymbol{\hat{\eta}}$ when $s$ is fixed and when  $s\to \infty$}\label{UpBnd}

The next result states that if the elements of the  family of  collections $r_\varepsilon=\{r_{\varepsilon,k},1\leq k\leq s\}$ in the hypothesis testing problem (\ref{HT2})
all exceed a certain level, then the upper bound on the maximum Hamming risk of $\boldsymbol{\hat{\eta}}$ as  defined in (\ref{r11}) holds true, and hence
the proposed selector $\boldsymbol{\hat{\eta}}$ is \textit{exact}.

	\begin{theorem} \label{Theorem1}
		Let $s\in\{1,\ldots,d\}$, $\beta\in(0,1)$, and $\sigma>0$ be fixed numbers, and let $d=d_\varepsilon\to \infty$
		and  $\log {d\choose s}=o\left( \log \varepsilon^{-1} \right)$ as $\varepsilon\to 0$.
		Let the family of collections  $r_\varepsilon=\{r_{\varepsilon,k},1\leq k\leq s\}$, $r_{\varepsilon,k}>0$, be such that
		\begin{equation}
			\liminf_{\varepsilon \to 0} \min_{1 \leq k \leq s} \frac{a_{\varepsilon,u_k}(r_{\varepsilon,k})}{\sqrt{\log {d\choose k}}} > \sqrt{2} (1+\sqrt{1-\beta}).
			\label{cond1}
		\end{equation}
		Then the selector $\boldsymbol{\hat{\eta}}=(\boldsymbol{\hat{\eta}}_k=(\hat{\eta}_{u_k},u_k \in {\cal U}_{k,d}),1\leq k\leq s)$ given by {\rm (\ref{selector})} and {\rm (\ref{delta})} satisfies
		\begin{equation*}
			\limsup_{\varepsilon\to 0}	\sup_{\boldsymbol{\eta} \in \mathcal{H}^s_{\beta,d}} \sup_{\thetab\in  \Theta_{s,d}^{\sigma}(r_{\varepsilon})} \E{\thetab,\etab} | \boldsymbol{\hat{\eta}}-\boldsymbol{\eta} | = 0.
		\end{equation*}
	\end{theorem}
	
\noindent	\textbf{Proof.}
	Denote  the threshold in (\ref{selector}) by $t_{\varepsilon,k}=\sqrt{(2+\epsilon) \left( \log {d\choose k} + \log M_k \right)},$ $1\leq k\leq s.$
When $\eta_{u_k}=0$ (respectively $\eta_{u_k}=1$), we shall write $\operatorname{P}_{0}$ and $\operatorname{E}_{0}$
(respectively $\operatorname{P}_{\thetab_{u_k}}$ and $\operatorname{E}_{\thetab_{u_k}}$) for $\operatorname{P}_{\thetab_{u_k},\eta_{u_k}}$ and $\operatorname{E}_{\thetab_{u_k},\eta_{u_k}}$.
For all small enough $\varepsilon$, the maximum Hamming risk $R_{\varepsilon,s}(\boldsymbol{\hat{\eta}})$ of the selector
$\boldsymbol{\hat{\eta}}=(\boldsymbol{\hat{\eta}}_k=(\hat{\eta}_{u_k},u_k \in {\cal U}_{k,d}),1\leq k\leq s)$
satisfies
\begin{align}
R_{\varepsilon,s}(\boldsymbol{\hat{\eta}})&=\sup_{\boldsymbol{\eta} \in  \mathcal{H}^s_{\beta,d}} \sup_{\thetab\in  \Theta_{s,d}^{\sigma}(r_{\varepsilon})} \E{\thetab,\etab} | \boldsymbol{\hat{\eta}}-\boldsymbol{\eta} | =\sup_{\boldsymbol{\eta} \in \mathcal{H}^s_{\beta,d}} \sup_{\thetab\in  \Theta_{s,d}^{\sigma}(r_{\varepsilon})}
\sum_{k=1}^s \sum_{u_k\in{\cal U}_{k,d}} \operatorname{E}_{\thetab_{u_k},\eta_{u_k}}|\hat{\eta}_{u_k}-{\eta}_{u_k}|\nonumber\\
&=\sup_{\boldsymbol{\eta} \in \mathcal{H}^s_{\beta,d}} \sup_{\thetab\in \Theta_{s,d}^{\sigma}(r_{\varepsilon})}
\sum_{k=1}^s \left(\sum_{u_k:\eta_{u_k}=0}  \operatorname{E}_{0}(\hat{\eta}_{u_k})+
\sum_{u_k:\eta_{u_k}=1} \operatorname{E}_{\thetab_{u_k}}(1-\hat{\eta}_{u_k})  \right) \nonumber\\
&\leq \sup_{\boldsymbol{\eta} \in \mathcal{H}^s_{\beta,d}}
\sum_{k=1}^s \sum_{u_k:\eta_{u_k}=0}  \operatorname{E}_0\left(\max_{1\leq m\leq M_k}\mathds{1}\left(S_{u_k,m}> t_{\varepsilon,k} \right)\right)\nonumber\\
&\quad \quad + \sup_{\boldsymbol{\eta} \in \mathcal{H}^s_{\beta,d}} \sup_{\thetab\in  \Theta_{s,d}^{\sigma}(r_{\varepsilon})}
\sum_{k=1}^s\sum_{u_k:\eta_{u_k}=1} \operatorname{E}_{\thetab_{u_k}}\left(1-\max_{1\leq m\leq M_k}\mathds{1}\left(S_{u_k,m}> t_{\varepsilon,k} \right)\right) \nonumber\\
&\leq \sum_{k=1}^s \sup_{\boldsymbol{\eta}_k \in  {H}^k_{\beta,d}}\sum_{u_k:\eta_{u_k}=0}\sum_{m=1}^{M_k}
\operatorname{P}_0\left(S_{u_k,m}> t_{\varepsilon,k} \right)\nonumber \\
& \quad \quad +
\sum_{k=1}^s \sup_{\boldsymbol{\eta}_k \in  {H}^k_{\beta,d}}\sum_{u_k:\eta_{u_k}=1}\sup_{\thetab_{u_k}\in \mathring{\Theta}_{\cb_{u_k}}(r_{\varepsilon,k})}
\operatorname{E}_{\thetab_{u_k}}\left\{\mathds{1}\left(\bigcap_{m=1}^{M_k}\{S_{u_k,m}\leq t_{\varepsilon,k}\}   \right)  \right\}\nonumber\\
&\leq \sum_{k=1}^s {d\choose k}\sum_{m=1}^{M_k}\operatorname{P}_0\left(S_{u_k,m}> t_{\varepsilon,k} \right)+
2\sum_{k=1}^s{d\choose k}^{1-\beta}\!\!\!\!\!\sup_{\thetab_{u_k}\in \mathring{\Theta}_{\cb_{u_k}}(r_{\varepsilon,k})}\min_{1\leq m\leq M_k}
\operatorname{P}_{\thetab_{u_k}}\left(S_{u_k,m}\leq t_{\varepsilon,k} \right)\label{upperbound0}.
\end{align}
By taking the index $m_{0,k}$ ($1\leq m_{0,k}\leq M_k$) to be such that $\beta\in(\beta_{m_{0},k}  ,\beta_{m_{0}+1,k}]$ for $1\leq k\leq~s$, we can further estimate the second term on the right-hand side of (\ref{upperbound0}) and arrive at
\begin{align}
R_{\varepsilon,s}(\boldsymbol{\hat{\eta}}) \leq& \sum_{k=1}^s {d\choose k}\sum_{m=1}^{M_k}\operatorname{P}_0\left(S_{u_k,m}> t_{\varepsilon,k} \right)+
2\sum_{k=1}^s{d\choose k}^{1-\beta}\!\!\!\!\!\sup_{\thetab_{u_k}\in \mathring{\Theta}_{\cb_{u_k}}(r_{\varepsilon,k})}
\operatorname{P}_{\thetab_{u_k}}\left(S_{u_k,m_{0,k}}\leq t_{\varepsilon,k} \right) \nonumber \\
=&: I^{(1)}_{\varepsilon,s}+I^{(2)}_{\varepsilon,s}.
\label{upperbound}
\end{align}
Similar to how this is done in the proof of Theorem 3.1 in Stepanova and Turcicova (2025), i.e., applying the exponential upper bound
(\ref{bound_1}) to the probability
$\operatorname{P}_0\left(S_{u_k,m}> t_{\varepsilon,k} \right)$, in view of (\ref{delta}) and the requirement $\log M_k=o(\log {d\choose k})$, we obtain as $\varepsilon\to 0$
\begin{eqnarray}
I^{(1)}_{\varepsilon,s}&\leq& \sum_{k=1}^s {d\choose k}\sum_{m=1}^{M_k} \exp\left( -(1+\epsilon/2)\left\{ \log {d\choose k}+\log M_k \right\}(1+o(1)) \right)\nonumber\\
&=&O\left(\sum_{k=1}^s \left\{M_k{d\choose k}  \right\}^{-\epsilon/2}  \right)=O\left(\max_{1\leq k\leq s}\left\{M_k{d\choose k}  \right\}^{-\epsilon/2}  \right)=o(1). \label{I1_eps}
\end{eqnarray}
Next, under the assumptions of the theorem, applying, as in the proof of Theorem 3.1 in  Stepanova and Turcicova (2025),
the exponential upper bound
(\ref{bound_2}) to the probability $\operatorname{P}_{\thetab_{u_k}}\left(S_{u_k,m_{0,k}}\leq t_{\varepsilon,k} \right)=
\operatorname{P}_{\thetab_{u_k}}\left(S_{u_k,m_{0,k}}-\operatorname{E}_{\thetab_{u_k}}(S_{u_k,m_{0,k}})\leq t_{\varepsilon,k}-\operatorname{E}_{\thetab_{u_k}}(S_{u_k,m_{0,k}}) \right)$
for those $\thetab_{u_k}\in \mathring{\Theta}_{\cb_{u_k}}(r_{\varepsilon,k})$ for which
$$\operatorname{E}_{\thetab_{u_k}}(S_{u_k, m_{0,k}}  )\max_{\lb\in \mathring{\mathbb{Z}}_{u_k}}\omega_{\lb}(r^*_{\varepsilon,k,m_{0,k}})
=o(1),\quad\varepsilon\to 0,$$ and Chebyshev's inequality for those $\thetab_{u_k}\in \mathring{\Theta}_{\cb_{u_k}}(r_{\varepsilon,k})$
for which the above limiting relation does not hold true, we obtain that
 $I^{(2)}_{\varepsilon,s}=o(1)$ as $\varepsilon\to 0$.
From this, (\ref{upperbound}), and (\ref{I1_eps}), $\limsup_{\varepsilon\to 0}R_{\varepsilon,s}(\boldsymbol{\hat{\eta}})=0$, which completes the proof.
\qedsymbol
	
\bigskip

	Theorem \ref{Theorem1} can be extended to the case of growing $s$.

	\begin{theorem} \label{Theorem2}
		Let  $\beta\in(0,1)$ and $\sigma>0$ be fixed numbers, and let
		$d=d_\varepsilon \to \infty$ and $s=s_\varepsilon \to \infty$ be such that
		$s = o(d)$, $\log \log d = o(s)$, and $\log {d\choose s}=o \left( \log \varepsilon^{-1} \right)$ as $\varepsilon\to 0$.
		Let the family of collections  $r_\varepsilon=\{r_{\varepsilon,k},1\leq k\leq s\}$, $r_{\varepsilon,k}>0$, be as in Theorem {\rm\ref{Theorem1}}.
		Then the selector $\boldsymbol{\hat{\eta}}=(\boldsymbol{\hat{\eta}}_k=(\hat{\eta}_{u_k},u_k \in {\cal U}_{k,d}),1\leq k\leq s)$ given by {\rm (\ref{selector})}
		and {\rm (\ref{delta1})} satisfies
		\begin{equation*}
			\limsup_{\varepsilon\to 0}	\sup_{\boldsymbol{\eta} \in  \mathcal{H}^s_{\beta,d}} \sup_{\thetab\in \Theta_{s,d}^{\sigma}(r_{\varepsilon})} \E{\thetab,\etab} | \boldsymbol{\hat{\eta}}-\boldsymbol{\eta}| = 0.
		\end{equation*}
	\end{theorem}
\noindent\textbf{Proof.}
The proof goes along the same lines as that of Theorem \ref{Theorem1}, with condition (\ref{delta1})
used instead of condition (\ref{delta}). For example, acting as in (\ref{I1_eps}), in view of (\ref{delta1})
and the requirement $\log M_k=o(\log {d\choose k})$, we obtain as $\varepsilon\to 0$
\begin{eqnarray*}
I^{(1)}_{\varepsilon,s}&\leq& \sum_{k=1}^s {d\choose k}\sum_{m=1}^{M_k} \exp\left( -(1+\epsilon/2)\left\{ \log {d\choose k}+\log M_k \right\}(1+o(1)) \right)\\
&=&O\left(\sum_{k=1}^s \left\{M_k{d\choose k}  \right\}^{-\epsilon/2}  \right)=O\left(s\max_{1\leq k\leq s}\left\{M_k{d\choose k}  \right\}^{-\epsilon/2}  \right)=o(1).
\end{eqnarray*}
Also, under the assumptions of Theorem \ref{Theorem2}, acting as in the proof of Theorem \ref{Theorem1}, we can obtain that  $I^{(2)}_{\varepsilon,s}=o(1)$ as $\varepsilon\to 0$, and hence
$\limsup_{\varepsilon\to 0}R_{\varepsilon,s}(\boldsymbol{\hat{\eta}})\leq
\limsup_{\varepsilon\to 0}\left(I^{(1)}_{\varepsilon,s}+I^{(2)}_{\varepsilon,s}\right)=0.$
The importance of condition $\log\log d=o(s)$ when $s=s_\varepsilon\to \infty$, $s=o(d)$ as $\varepsilon\to 0$ can be seen from the proof
of Theorem~3.3 in Stepanova and Turcicova (2025). \qedsymbol

\bigskip

Theorems  \ref{Theorem1} and \ref{Theorem2}  show that
	if a positive family of collections $r_\varepsilon=\{r_{\varepsilon,k}:1\leq k\leq s\}$ is such that
  condition (\ref{cond1}) holds, exact variable selection in the sequence space model (\ref{model2_s}) is possible;
 this is true for fixed $s$ as well as for $s=s_\varepsilon\to \infty$ $s=o(d)$ as $\varepsilon\to 0$, provided the rate at which $s$ tends to infinity
  is carefully regulated.
	Theorems \ref{Theorem1} and \ref{Theorem2}  may be viewed as extensions of  Theorems 3.1 and 3.3 in  Stepanova and Turcicova (2025), for more details we refer to
the discussion of Section \ref{Discussion}.

\subsection{Lower bound on the minimax risk when $s$ is fixed and when $s\to \infty$}

We now turn to deriving the conditions when exact variable selection in
model (\ref{model2_s}) is impossible, i.e. when the lower bound (\ref{r22}) holds true, which occurs if the family of collections
$r_\varepsilon=\{r_{\varepsilon,k},1\leq k\leq s\}$ in the hypothesis testing problem (\ref{HT2}) falls below a certain level.
The precise statement is as follows.

\begin{theorem} \label{Theorem3}
	Let $s\in\{1,\ldots,d\}$, $\beta\in(0,1)$, and $\sigma>0$ be fixed numbers, and let $d=d_\varepsilon\to \infty$
	and $\log {d\choose s}=o \left(\log \varepsilon^{-1} \right)$ as $\varepsilon\to 0$.
	Let the family of collections  $r_\varepsilon=\{r_{\varepsilon,k},1\leq k\leq s\}$, $r_{\varepsilon,k}>0$, be such that
	\begin{equation}
		\limsup_{\varepsilon \to 0} \min_{1 \leq k \leq s} \frac{a_{\varepsilon,u_k}(r_{\varepsilon,k})}{\sqrt{\log {d\choose k}}} < \sqrt{2} (1+\sqrt{1-\beta}).
		\label{cond:inf2}
	\end{equation}
	Then
	\begin{equation*}
		\liminf_{\varepsilon \to 0}\inf_{\tilde{\etab}} \sup_{\boldsymbol{\eta} \in \mathcal{H}^s_{\beta,d}} \sup_{\thetab\in \Theta_{s,d}^{\sigma}(r_{\varepsilon})} \E{\thetab,\etab} | \boldsymbol{\tilde{\eta}}-\boldsymbol{\eta} |>0,
	\end{equation*}
	where the infimum is taken over all selectors $\tilde{\etab}$ of $\etab\in  \mathcal{H}^s_{\beta,d}$
	in model {\rm (\ref{model2_s})}.
\end{theorem}

\noindent \textbf{Proof.}
For every $k$ ($1\leq k\leq s$), we choose some $u_k\in {\cal U}_{k,d}$ and fix it.
Let $k^{\prime}=k^{\prime}(\varepsilon)$ be a map from $(0,\infty)$ to $\{1,\ldots,s\}$ defined as follows:
\begin{equation*}
	k^{\prime}=\arg\!\min_{\!\!\!\!\!\!\!\!\!\!1\leq k\leq s}\frac{a_{\varepsilon,u_k}(r_{\varepsilon,k})}{\sqrt{\log {d\choose k}}}.
\end{equation*}
We have
\begin{align}
	R_{\varepsilon,s}:=&\inf_{\tilde{\etab}} \sup_{\boldsymbol{\eta} \in  \mathcal{H}^s_{\beta,d}} \sup_{\thetab\in \Theta_{s,d}^{\sigma}(r_{\varepsilon})} \E{\thetab,\etab} | \boldsymbol{\tilde{\eta}}-\boldsymbol{\eta}|\nonumber \\
	=&\inf_{\tilde{\etab}} \sup_{\boldsymbol{\eta} \in  \mathcal{H}^s_{\beta,d}} \sup_{\thetab\in\Theta_{s,d}^{\sigma}(r_{\varepsilon})}  \E{\thetab,\etab} \left(\sum_{k=1}^{s}\sum_{u_k \in \mathcal{U}_{k,d}} | \tilde{\eta}_{u_k}-{\eta}_{u_k}|\right)\nonumber\\
	\geq&
	\inf_{\tilde{\etab}_{k^{\prime}}} \sup_{\boldsymbol{\eta}_{k^{\prime}} \in {H}^{k^{\prime}}_{\beta,d}} \sup_{\thetab_{k^{\prime}}\in \mathring{\Theta}^{\sigma}_{k^{\prime},d}(r_{\varepsilon,k^{\prime}})}
	\E{\thetab_{k^{\prime}},\etab_{k^{\prime}}} \left(\sum_{u_{k^{\prime}} \in \mathcal{U}_{k^{\prime},d}} | \tilde{\eta}_{u_{k^{\prime}}}-{\eta}_{u_{k^{\prime}}}|\right).\label{eq1}
\end{align}
Since $\log {d\choose s}=o \left(\log \varepsilon^{-1} \right)$
implies $\log {d\choose k^{\prime}}=o \left(\log \varepsilon^{-1} \right)$, it follows from (\ref{eq1}) and Theorem 3.2 in Stepanova and Turcicova (2025)
that
$$\liminf_{\varepsilon\to 0}R_{\varepsilon,s}>0$$
provided
\begin{gather*}
	\limsup_{\varepsilon \to 0} \frac{a_{\varepsilon,u_{k^{\prime}}}(r_{\varepsilon,k^{\prime}})}{\sqrt{\log {d\choose k^{\prime}}}} < \sqrt{2} (1+\sqrt{1-\beta}),
\end{gather*}
which is true by the definition of $k^{\prime}$ and condition (\ref{cond:inf2}). This completes the proof. \qedsymbol
\bigskip

The result of Theorem \ref{Theorem3} can be extended to the case of growing $s$.

\begin{theorem} \label{Theorem4}
	Let  $\beta\in(0,1)$ and $\sigma>0$ be fixed numbers,
	and let  $d=d_\varepsilon \to \infty$ and $s=s_\varepsilon \to \infty$ be such that
	$s = o(d)$, $\log \log d = o(s)$, and $\log {d\choose s}=o \left( \log \varepsilon^{-1} \right)$ as $\varepsilon\to 0$.
	Let the family of collections  $r_\varepsilon=\{r_{\varepsilon,k},1\leq k\leq s\}$, $r_{\varepsilon,k}>0$, be as in Theorem {\rm \ref{Theorem3}}.
	Then
	\begin{equation*}
		\liminf_{\varepsilon \to 0}\inf_{\tilde{\etab}} \sup_{\boldsymbol{\eta} \in  \mathcal{H}^s_{\beta,d}} \sup_{\thetab\in \Theta_{s,d}^{\sigma}(r_{\varepsilon})} \E{\thetab,\etab} | \boldsymbol{\tilde{\eta}}-\boldsymbol{\eta}|>0,
	\end{equation*}
	where the infimum is taken over all selectors $\tilde{\etab}$ of $\etab\in  \mathcal{H}^{s}_{\beta,d}$
	in model {\rm (\ref{model2_s})}.
\end{theorem}

\noindent\textbf{Proof.}
Using the same arguments as in the proof of Theorem \ref{Theorem3},
with reference to Theorem~3.4 of Stepanova and Turcicova (2025) instead of reference to Theorem 3.2 of  Stepanova and Turcicova (2025), we arrive at the statement of Theorem~\ref{Theorem4}.  \qedsymbol
\bigskip

Theorems \ref{Theorem3} and \ref{Theorem4} show that
	if a positive family of collections $r_\varepsilon=\{r_{\varepsilon,k}:1\leq k\leq s\}$ is such that
  condition (\ref{cond:inf2}) holds, exact variable selection in the sequence space model (\ref{model2_s}) is impossible;
 this is true for fixed $s$ as well as for $s=s_\varepsilon\to \infty$ $s=o(d)$ as $\varepsilon\to 0$, provided the rate at which $s$ tends to infinity
  is carefully regulated.
Theorems \ref{Theorem3} and \ref{Theorem4}  may be viewed as extensions of  Theorems 3.2 and~3.4 in  Stepanova and Turcicova (2025),  for more details
we refer to the discussion in Section \ref{Discussion}.

\section{Discussion}\label{Discussion}
By using the notion of optimality of a statistical procedure from the minimax hypothesis testing theory, we arrive at the following important conclusions.
When $s$ is fixed, in view of Theorems \ref{Theorem1} and~\ref{Theorem3},   the selector
defined in {\rm (\ref{selector})} and~{\rm (\ref{delta})} is seen to perform
\textit{optimally} with respect to a Hamming loss in the asymptotically minimax sense.
Similarly, when $s=s_\varepsilon\to \infty$, $s=o(d)$ as $\varepsilon\to 0$, by means of Theorems \ref{Theorem2} and \ref{Theorem4}, the selector
defined in {\rm (\ref{selector})} and {\rm (\ref{delta1})} is optimal with respect to a Hamming loss in the {asymptotically minimax} sense.

We now explain our choice to state the sharp selection boundary in terms of the
 quantities $a_{\varepsilon,u_k}(r_{\varepsilon,k})$, $1\leq k\leq s$, rather than in terms of
the radii $r_{\varepsilon,k}$, $1\leq k\leq s$. To this end, consider the Gaussian vector model
\begin{gather}\label{vm}
	X_{u_k}=\mu_{k,d}\eta_{u_k}+\varepsilon_{u_k},\quad u_k\in {\cal U}_{k,d},\quad 1\leq k\leq s,
\end{gather}
where  $\mu_{k,d}>0$ is the signal, the errors $\varepsilon_{u_k}$ are iid standard normal random variables, the deterministic quantities $\eta_{u_k}\in \{0,1\}$ are as before and, in particular, satisfy the sparsity condition (\ref{SC1}).
 For $\eta_{u_k}=1$, the observation $X_{u_k}$ is equal
to the signal $\mu_{k,d}$ perturbed by random noise $\varepsilon_{u_k}$.
For  $\eta_{u_k}=0$, the observation $X_{u_k}$ is merely random noise $\varepsilon_{u_k}$.
Due to condition (\ref{SC1}), only
$\sum_{k=1}^{s} \sum_{u_k \in \mathcal{U}_{k,d}} \eta_{u_k}= \sum_{k=1}^s {d\choose k}^{1-\beta}(1+o(1))= {d\choose s}^{1-\beta}(1+o(1))$
observations in model ({\ref{vm}) contain signal, and thus the model is sparse.
The sharp selection boundary in model (\ref{vm}) is given by, cf. inequalities (6) and (7) in Stepanova and Turcicova (2025),
\begin{gather*}
\liminf_{d \to \infty} \min_{1\leq k\leq s} \frac{\mu_{k,d}}{\sqrt{\log {d\choose k}}} > \sqrt{2} (1+\sqrt{1-\beta})\quad\mbox{and}\quad
\limsup_{d \to \infty} \min_{1\leq k\leq s}\frac{\mu_{k,d}}{\sqrt{\log {d\choose k}}} < \sqrt{2} (1+\sqrt{1-\beta}).
\end{gather*}
By comparing the above inequalities to the sharp selection boundary determined by (\ref{cond1}) and  (\ref{cond:inf2}), we conclude that,
 in the variable selection problem at hand, the quantity $a_{\varepsilon,u_k}(r_{\varepsilon,k})$
plays the same role as the signal $\mu_{k,d}$ does in the problem of recovering
$\boldsymbol{{\eta}}=(\boldsymbol{{\eta}}_k=({\eta}_{u_k},u_k \in {\cal U}_{k,d}),1\leq k\leq s)$ in model~(\ref{vm}).
Additionally, the use of $a_{\varepsilon,u_k}(r_{\varepsilon,k})$, $1\leq k\leq s$, for establishing the sharp selection boundary
is natural for comparison purposes, as seen from the discussion below.

Let us compare the main results of this article to those stated
 in Theorems 3.1 to 3.4  of  Stepanova and Turcicova (2025) for
a simpler model when, instead of decomposition (\ref{fun_sub_anova22}) for function $f$ in model (\ref{model1})--(\ref{SC1}) and (\ref{Fu}), one has decomposition (\ref{f}).
The \textit{sharp selection boundary} for the latter (simpler) model is stated in terms of  $a_{\varepsilon,u_k}(r_{\varepsilon,k})$ as in (\ref{def:a2}) and is given by the inequalities,
cf. inequalities {\rm (\ref{cond1})} and  {\rm (\ref{cond:inf2})},
\begin{gather*}
			\liminf_{\varepsilon \to 0}  \frac{a_{\varepsilon,u_k}(r_{\varepsilon,k})}{\sqrt{\log {d\choose k}}} > \sqrt{2} (1+\sqrt{1-\beta})\quad
\mbox{and}\quad\limsup_{\varepsilon \to 0}  \frac{a_{\varepsilon,u_k}(r_{\varepsilon,k})}{\sqrt{\log {d\choose k}}} < \sqrt{2} (1+\sqrt{1-\beta}).
\end{gather*}
Here $u_k$ is an arbitrary element of ${\cal U}_{k,d}$, and $r_{\varepsilon,k}$ determines the set $\mathring{\Theta}_{\cb_{u_k}}(r_{\varepsilon,k})$ as defined in (\ref{Theta}),
more precisely, $r_{\varepsilon,k}$ is a radius of the centered $\l2$-ball removed from the ellipsoid ${\Theta}_{\cb_{u_k}}$.
The above inequalities provide conditions on the $\l2$-norm of $(\theta_{\lb}(u_k))_{\lb \in \mathring{\mathbb{Z}}_{u_k}}$, $u_k\in {\cal U}_{k,d}$,
 that make exact variable selection in the Gaussian sequence space model studied in Stepanova and Turcicova (2025) possible (the first inequality) and impossible
 (the second inequality).
In terms of the initial model (\ref{model1})--(\ref{SC1}) and (\ref{Fu}), with (\ref{f}) in place of (\ref{fun_sub_anova22}),
the quantity $r_{\varepsilon,k}$ determines the set $\mathring{\cal F}_{\cb_{u_k}}(r_{\varepsilon,k})=\{f_{u_k}\in {\cal F}_{\cb_{u_k}}: \|f_{u_k}\|_2\geq r_{\varepsilon,k}\} $,
and thus the above inequalities also provide conditions on the
$L_2$-norm of the components $f_{u_k}$ of $f=\sum_{u_k\in {\cal U}_{k,d}} \eta_{u_k}f_{u_k}$ that make exact selection of the active (nonzero) components possible
and impossible.
Note that the sharp selection boundary is the same for fixed $k$ and for $k\to \infty$, $k=o(d)$ as $\varepsilon\to 0$
provided the rate at which $k$ tends to infinity is carefully regulated (see Theorems 3.1 to 3.4 in Stepanova and Turcicova (2025)).
This is consistent with the results in Section 3 of this article when  $s$ is fixed (see Theorems \ref{Theorem1} and \ref{Theorem3}) and when $s\to \infty$, $s=o(d)$ as $\varepsilon\to 0$
(see Theorems \ref{Theorem2} and \ref{Theorem4}).
In case of growing $s$, the adaptive selector $\boldsymbol{\hat{\eta}}=(\boldsymbol{\hat{\eta}}_k=(\hat{\eta}_{u_k},u_k \in {\cal U}_{k,d}),1\leq k\leq s)$,
where  ${\hat{\eta}}_{u_k}$ is defined by~(\ref{selector}),
is similar to that for fixed $s$, with a slight difference in the threshold $\sqrt{(2+\epsilon) \left( \log {d\choose k} + \log M_k \right)}$,
as seen from conditions (\ref{delta}) and (\ref{delta1}) imposed on $\epsilon$.

It is also of interest to make a bridge between the sharp selection boundary given by inequalities
{\rm (\ref{cond1})} and  {\rm (\ref{cond:inf2})} and the \textit{sharp detection boundary} established by Ingster and Suslina (2015)
in the problem of testing
\begin{gather*}
\mathbf{H}_{0,s}:f=0\quad\mbox{vs.}\quad \mathbf{H}_{1,s}^{\varepsilon}: f\in\bigcup\limits_{1\leq k\leq s}\bigcup\limits_{u_k\in{\cal U}_{k,d}}\mathring{\cal F}_{\cb_{u_k}}(r_{\varepsilon,k}).
\end{gather*}
As shown in Theorem 3.3 of Ingster and Suslina (2015), the sharp detection boundary is given by the inequalities
\begin{equation}\label{testing5}
	\liminf_{\varepsilon \to 0} \min_{1\leq k\leq s} \frac{a_{\varepsilon,u_k}(r_{\varepsilon,k})}{\sqrt{\log {d\choose k}}} > \sqrt{2}
	\quad \mbox{and}\quad  \limsup_{\varepsilon \to 0} \min_{1\leq k\leq s} \frac{a_{\varepsilon,u_k}(r_{\varepsilon,k})}{\sqrt{\log {d\choose k}}} < \sqrt{2}.
\end{equation}
Specifically, when the first inequality in (\ref{testing5}) holds,
the hypotheses $\mathbf{H}_{0,s}$ and $ \mathbf{H}_{1,s}^{\varepsilon}$ separate asymptotically (i.e., there exists a consistent test
procedure for testing $\mathbf{H}_{0,s}$ versus $ \mathbf{H}_{1,s}^{\varepsilon}$), whereas, when the second inequality in (\ref{testing5}) holds,
$\mathbf{H}_{0,s}$ and $ \mathbf{H}_{1,s}^{\varepsilon}$ merge asymptotically (i.e., there is no consistent test procedure
for testing $\mathbf{H}_{0,s}$ versus $ \mathbf{H}_{1,s}^{\varepsilon}$).
By comparing the inequalities in (\ref{testing5}) to the inequalities {\rm (\ref{cond1})} and  {\rm (\ref{cond:inf2})},
we notice that the sharp selection boundary lies above the sharp detection boundary. This is not unusual because
the problem of selecting nonzero components is more difficult than that of testing whether such nonzero components exit.
In particular, in order to be \textit{selectable}, the components need to be \textit{detectable}.

\section{Simulation study}
To showcase the exact selection procedure {based on} $\boldsymbol{\hat{\eta}}=(\boldsymbol{\hat{\eta}}_k=(\hat{\eta}_{u_k},u_k \in {\cal U}_{k,d}),1\leq k\leq s)$ in (\ref{selector})--(\ref{delta}) in its capacity
to recover the sparsity pattern of a signal $f\in {\cal F}^{\beta,\sigma}_{s,d}(r_{\varepsilon})$
of the form $f=\sum_{k=1}^s \sum_{u_k\in {\cal U}_{k,d}}\eta_{u_k} f_{u_k}$, we conduct a small-scale simulation study.
For this study, we choose $\sigma=1$, $s=4$, and $d=50, 100, 200$. We also choose $\beta = 0.87$, which
means that the signal $f$ is highly sparse.
The number of active components $ N_{k,d} := \sum_{u_k \in {\cal U}_{k,d} } \eta_{u_k}$ of $f$ for each $k$ is obtained from relation (\ref{SC1}). The  values of $N_{k,d}$ are listed in Table \ref{tab:number_of_actives}.

\begin{table}[h]
	\centering
	\caption{The number $N_{k,d}$ of active components of $f$ for $\beta=0.87$ and the selected values of $d$ and $k$.}
	\label{tab:number_of_actives}
	\begin{tabular}{@{}lc|c|c|c|c@{}}
		\hline
		& &$k=1$ & $k=2$ &$k=3$ & $k=4$\\			
		\hline
		\multirow{3}{*}{$d$} & 50 & 2 & 3 &4 &5 \\
		&100 & 2 & 3 &5&7\\
		&200 & 2 & 3 &6&10\\
		\hline
	\end{tabular}
\end{table}

In this simulation study, we choose the active components of $f$ as follows. First, consider the following nine functions defined on $[0,1]$:
\begin{align*}
	g_1(t) &= t^2 (2^{t-1} - (t-0.5)^2 ) \exp(t) - 0.5424 , \\
	g_2(t) &= t^2 (2^{t-1} - (t-1)^5) - 0.2887, \\
	g_3 (t)&= 1.5 t^2 2^{t-1} \cos(15 t) - 0.05011, \\
	g_4 (t)&= t-{1}/{2}, \\
	g_5 (t)&= 5 (t-0.7)^3 + 0.29,\\
	g_6(t)&=2(t-0.4)^2-0.1867,\\
	g_7(t)&=0.7(t^2-0.1)^3-0.0643,\\
	g_8(t)&=10(t^2-0.5)^5 + 0.068,\\
	g_9(t)&=3(t-0.8)^4-0.1968.
\end{align*}
In order to simplify the notation in this section, instead of indexing the components $f_{u_k}$ of $f$ by $u_k$, (e.g., $f_{\{i,j\}}$ for $u_2=\{i,j\} $),
we will index them by   $u_k^{(i)}$ for $i = 1,\ldots, {d \choose k}$.
By using this notation, the active (and inactive) components $f_{u_k}$ of $f$  on $[0,1]^k$, $k=1,\ldots,4$, are defined as follows:
\begin{align*}
	\begin{split}
		\text{for $k=1$:} \quad &f_{u_{1}^{(i)}}(t_i)=g_{i}(t_i), \quad i=1,2, \\
		&f_{u_1^{(i)}}(t_{u_1^{(i)}})=0, \quad i = 3,\ldots, d, \\
		\text{for $k=2$:} \quad &f_{u_{2}^{(i)}}(t_{i},t_{i+1})=g_{i}(t_i)g_{i+1}(t_{i+1}), \quad i=1,2,3, \\
		&f_{u_2^{(i)}}(\tb_{u_2^{(i)}})=0, \quad i =4,\ldots, {d \choose 2}, \\
		\text{for $k=3$:} \quad & f_{u_{3}^{(i)}}(t_1,t_2,t_{{i+2}})=g_1(t_1)g_2(t_2)g_{i+2}(t_{i+2}), \quad i=1,\ldots,4, \\
		& \text{for $d = 50$: } f_{u_{3}^{(i)}}(\tb_{u_3^{(i)}})=0, \quad i = 5,\ldots, {d \choose 3},\\
		&\text{for $d = 100$: } f_{u_{3}^{(5)}}(t_1,t_2,t_7)=g_1(t_1)g_2(t_2)g_{7}(t_7), \\
		& \hspace{2.1cm} f_{u_{3}^{(i)}}(\tb_{u_3^{(i)}})=0, \quad i = 6,\ldots, {d \choose 3},\\
		&\text{for $d = 200$: } f_{u_{3}^{(i)}}(t_1,t_2,t_{i+2})=g_1(t_1)g_2(t_2)g_{i+2}(t_{i+2}), \quad i=5,6, \\
		& \hspace{2.1cm} f_{u_{3}^{(i)}}(\tb_{u_3^{(i)}})=0, \quad i = 7,\ldots, {d \choose 3},\\
		\text{for $k=4$:} \quad &f_{u_{4}^{(i)}}(t_1,t_2,t_3,t_{i+3})=g_1(t_1)g_2(t_2)g_3(t_3)g_{i+3}(t_{i+3}), \quad i=1,\ldots,5, \\
		& \text{for $d = 50$: } f_{u_{4}^{(i)}}(\tb_{u_4^{(i)}})=0, \quad i = 6,\ldots, {d \choose 4},\\
		&\text{for $d = 100$: } f_{u_{4}^{(i)}}(t_1,t_2,t_4,t_{i+2})=g_1(t_1)g_2(t_2)g_4(t_4)g_{i+2}(t_{i+2}), \quad i=6,7,\\
		& \hspace{2.1cm} f_{u_{4}^{(i)}}(\tb_{u_4^{(i)}})=0, \quad i = 8,\ldots, {d \choose 4},\\
		&\text{for $d = 200$: }  f_{u_{4}^{(i)}}(t_1,t_2,t_4,t_{i+2})=g_1(t_1)g_2(t_2)g_4(t_4)g_{i+2}(t_{i+2}), \quad i=6,7,\\
		& \hspace{2.1cm} f_{u_{4}^{(i)}}(t_1,t_2,t_5,t_{i-1})=g_1(t_1)g_2(t_2)g_5(t_5)g_{i-1}(t_{i-1}), \quad i=8,9,10,\\
		& \hspace{2.1cm} f_{u_{4}^{(i)}}(\tb_{u_4^{(i)}})=0, \quad i = 11, \ldots, {d \choose 4},\\
	\end{split}
	\label{def:functions_fu}
\end{align*}
where $\tb_{u_k^{(i)}} = (t_j)_{j \in u_k^{(i)}}$ as before.
When the component functions $f_{u_k^{(i)}}$ are as above,
the orthogonality condition (\ref{orthcon2}) holds true up to four decimal places
for all chosen values of $k$ and $d$.
The active components are associated with $u_k^{(i)} = \{i\}, \, i=1, 2,$ for $k=1$, $u_k^{(i)} = \{i,i+1\}, \, i= 1,2,3,$ for $k=2$, $u_k^{(i)} = \{1,2,i+2\}, \, i=1,\ldots,N_{3,d},$ for $k=3$ and $d=50, 100, 200$. For $k=4$, they are associated with  $u_k^{(i)}= \{1,2,3,i+3\}, \, i=1, \ldots, 5,$ when $d=50, 100, 200; $
$ u_k^{(i)} = \{1,2,4,i+2\}, \, i=6,7,$ when $d=100, 200;$ and $ u_k^{(i)} = \{1,2,5,i-1\}, \, i=8, 9, 10,$ when $d=200$.

Recall the definition of the exact selector $\boldsymbol{\hat{\eta}}=(\boldsymbol{\hat{\eta}}_k=(\hat{\eta}_{u_k},u_k \in {\cal U}_{k,d}),1\leq k\leq s)$ with components
$\hat{\eta}_{u_k}$ given by  (\ref{selector})--(\ref{delta}):
\begin{equation*}
	\hat{\eta}_{u_k} = \max_{1 \leq m \leq M_k}  \ind{ S_{u_k,m} > \sqrt{(2+\epsilon) \left( \log {d\choose k} + \log M_k \right)} },
\end{equation*}
where the weights of the statistics
$S_{u_k,m} = \sum_{\lb \in \mathring{\mathbb{Z}}_{u_k}} \omega_\lb (r^*_{\varepsilon,k,m}) \left( \left( {X_\lb}/{\varepsilon}  \right)^2 -1 \right)$,
$1\leq m\leq M_k$,
are equal to $$\omega_\lb (r^*_{\varepsilon,k,m}) = \dfrac{1}{2 \varepsilon^2} \dfrac{(\theta^*_\lb (r^*_{\varepsilon,k,m}))^2}{a_{\varepsilon,u_k}(r^*_{\varepsilon,k,m})},
\quad \lb\in  \mathring{\mathbb{Z}}_{u_k},$$
and the values of $\theta^*_\lb (r_{\varepsilon,k})$ and $a_{\varepsilon,u_k}(r_{\varepsilon,k})$ are computed by using relations (\ref{thetal}) and (\ref{aek}).
Keeping in mind condition (\ref{delta}) and the requirements on $M_k$, we take $\epsilon=1/\sqrt{\log {d \choose k}}$, $M_k=20$, and choose
the points $\beta_{k,m}$, $m=1, \ldots, M_k$, to be uniformly placed on the interval $[0.001, 0.999]$.
For each $k$ and $m$, the value $r^*_{\varepsilon,k,m}$ is found as a solution of equation (\ref{def:r_star_ekm}) with
$a_{\varepsilon,u_k} (r_{\varepsilon,k})$ as in (\ref{aek}).

In order to simulate the observations in model (\ref{model2_s}), we need to compute the true Fourier coefficients  $\theta_\lb(u_k), \lb \in \mathring{\mathbb{Z}}_{u_k}$,
in decomposition (\ref{f-decomp}). In what follows, $u_k$ stands for $u_k^{(i)}$ for $i =1,\ldots,N_{k,d}$.
When $k=1$, we simply have $u_1 = \{j\}$ for $j = 1,2$, $\tb_{u_1} = t_{j}$, $\lb = l_j$ and $$\theta_\lb\left(u_1 \right) = \left( f_{u_1}, \phi_{l_{j}}  \right)_{L_2^1} = ( g_{j} ,\phi_{l_{j}} )_{L^1_2}.$$
When $k=2$, we have $u_2=\{i,i+1\} =: \{j_1, j_2\}$, for $i=1,2,3$, $\tb_{u_k} = (t_{j_1},t_{j_2})$, and each index $\lb\in \mathring{\mathbb{Z}}_{u_k}$ has only two nonzero elements. Denote these nonzero elements of $\lb$ by $l_{j_1}$ and~$l_{j_2}$. Then
\begin{align*}
	\theta_\lb(u_2) &= ( f_{u_2}, \phi_{l_{j_1}}  \phi_{l_{j_2}} )_{L_2^2}  = \int_{0}^{1} \int_{0}^{1} g_{j_1}(t_{j_1}) g_{j_2} (t_{j_2})  \phi_{l_{j_1}} (t_{j_1}) \phi_{l_{j_2}} (t_{j_2}) dt_{j_1} dt_{j_2} \nonumber \\
	&= \int_{0}^{1}  g_{j_1}(t_{j_1})   \phi_{l_{j_1}} (t_{j_1})  dt_{j_1} \int_{0}^{1}g_{j_2} (t_{j_2}) \phi_{l_{j_2}} (t_{j_2})dt_{j_2}  = ( g_{j_1} ,\phi_{l_{j_1}} )_{L^1_2}  (g_{j_2},  \phi_{l_{j_2}} )_{L^1_2},
\end{align*}
and similarly for $k= 3$ and 4. For $k=1,2,3,4$, the Fourier coefficients $\theta_\lb(u_k^{(i)})$ with $i > N_{k,d}$ are all equal to zero.
For smooth Sobolev functions, the absolute values of their Fourier coefficients decay to zero at a~polynomial rate.
Therefore, although in theory $\lb=(l_1,\ldots,l_d) \in \mathring{\mathbb{Z}}_{u_k}$, we shall restrict ourselves to  $\lb \in \Ztuk \deq  \mathring{\mathbb{Z}}_{u_k} \cap [-n,n]^d,$ where $n=622$ for $k=1$, $n=154$ for $k=2$, $n=65$ for $k=3$, and $n=36$ for $k=4.$ For each $k$, the chosen value of $n$ ensures that none of the nonzero coefficients $\theta^*_\lb (r_{\varepsilon,k,m}^*)$, and hence none of the weights $\omega_\lb (r^*_{\varepsilon,k,m})$,  is missing in the evaluation of the statistics $S_{u_k,m}$, $m=1,\ldots, M_k$.
Thus, the model we are dealing with in this section  is as follows:
\begin{equation*}
	X_{\lb}=\eta_{u_k}\theta_\lb(u_k)+\varepsilon \xi_{\lb},\quad \lb \in \Ztuk, \quad 1 \leq k \leq 4,
\end{equation*}
with $\varepsilon=5 \cdot 10^{-5}$. The vector $(\xi_\lb)_{\lb \in \Ztuk}$ consists of iid standard normal random variables $\xi_{\lb}$,  and the component $\eta_{u_k} $ 
of $\boldsymbol{\eta}=(\boldsymbol{\eta}_k=(\eta_{u_k},u_k \in {\cal U}_{k,d}),1\leq k\leq 4)$ equals 1 if $u_k \in \left\{ u_k^{(1)},u_k^{(2)}, \ldots, u_k^{(N_{k,d})}  \right\}$  and zero otherwise.

For all choices of $d$, we run $J=15$ independent cycles of simulations (using R software) and
estimate the Hamming risk ${\rm E}_{\thetab,\etab}\left(\sum_{k=1}^s \sum_{u_k\in{\cal U}_{k,d}}|\hat{\eta}_{u_k}-\eta_{u_k}|\right)$
by means of the quantity
$${Err}(\hat{\etab})=\frac{1}{J}\sum_{j=1}^J\sum_{k=1}^s\sum_{u_k \in \mathcal{U}_{k,d}}|\hat{\eta}_{u_k}^{(j)}-\eta_{u_k}|,$$
where $\hat{\eta}_{u_k}^{(j)}$ is the value of $\hat{\eta}_{u_k}$ obtained in the $j$-th repetition of the experiment.
The values of ${Err}(\boldsymbol{\hat{\eta}})$ for different $d$ are seen in Table \ref{tab:Hast_err} in the last column corresponding to $\alpha=1$,
and they are all zero.

To illustrate the impact of a signal strength on the Hamming risk of $\boldsymbol{\hat{\eta}}$, we multiply the active component $f_{u_{1}^{(1)}}$ by $\alpha \in (0,1] $,
while keeping the other active components unchanged.
The values of the estimated Hamming risk ${Err}(\boldsymbol{\hat{\eta}})$ obtained for different choices of $\alpha$ are presented in Table~\ref{tab:Hast_err}.

\begin{table}[h!]
	\centering
	\caption{Estimated Hamming risk  ${Err}(\boldsymbol{\hat{\eta}})$ obtained from $J=15$ simulation cycles.}
	\label{tab:Hast_err}
	\begin{tabular}{@{}c|ccccccccc@{}}
		\hline
		& \multicolumn{9}{| c }{$\alpha$} \\
		 $d$& 0.0001 & 0.0005 & 0.0009 & 0.001 &  0.0011 & 0.0012 & 0.005 & 0.5 & 1  \\
		\hline
	  50 & 1 & 1 & 0.86 &  0.80 & 0.53 & 0.20  & 0  &0 &0 \\
	  100 &1 & 1 & 0.93 &  0.93 & 0.67 &  0.47 & 0  & 0 & 0\\
	  200 &1 & 1 & 1 &  1 & 0.93 &  0.67 & 0  &0 & 0\\
		\hline
	\end{tabular}
\end{table}

The selection procedure never detects a signal if there is none, and thus only active components may be misidentified as inactive.
Since, out of all the active components, only one of them is multiplied by $\alpha\in(0,1]$, the entries in Table \ref{tab:Hast_err}
(i.e., the average numbers of misidentified components of the signal $f=\sum_{k=1}^4 \sum_{u_k \in {\cal U}_{k,d}} \eta_{u_k} f_{u_k}$) are all less than or equal to 1.
In general, as seen from Table \ref{tab:Hast_err}, the stronger the signal is, the smaller is the number of misidentified active components.
It is also seen from  Table~\ref{tab:Hast_err} that the estimated Hamming risk ${Err}(\hat{\etab})$  increases as $d$ increases.
This is not surprising since the model gets sparser as $d$ gets larger,
and the sparser the model is, the harder is to select correctly the active components of $f$.
Overall, the numerical results of this section are  in agreement with the statement of Theorem \ref{Theorem1}.

\section{Concluding remarks}\label{CR}
The results obtained in this article provide the conditions for the possibility and impossibility of exact variable selection
with respect to a Hamming loss
in the Gaussian white noise model of intensity $\varepsilon$ when the $d$-variate signal $f$
is the sum of a small number of $k$-variate functions with $k$ varying from 1 to $s$ ($1\leq s\leq d$), for $s$ being fixed and for $s=s_\varepsilon\to \infty$, $s=o(d)$
as $\varepsilon\to 0$.
The {sharp selection boundary}, which defines a precise demarcation between what is possible and impossible in this problem, is determined by inequalities {\rm (\ref{cond1})} and  {\rm (\ref{cond:inf2})} in terms of the quantities $a_{\varepsilon,u_k}(r_{\varepsilon,k})$, $1\leq k\leq s$, defined by (\ref{def:a2})
and satisfying
$a_{\varepsilon,u_k}(r_{\varepsilon,k})\sim c(\sigma,k) r_{\varepsilon,k}^{2+k/(2\sigma)} \varepsilon^{-2}$ as $\varepsilon\to 0,$ where $c(\sigma,k)$ is a
known function of  $\sigma$ and $k$ (see relations (\ref{aek}) and (\ref{aek1})).
Under the conditions when exact variable selection is possible, we proposed an adaptive
selection procedure that achieves exact selection when the level of sparsity is unknown.

In this article, the regularity constraints on the $d$-variate function $f$ are expressed in terms of
the Sobolev semi-norms for $\sigma>0$ and are imposed component-wise. In addition to the Sobolev balls ${\cal F}_{\cb_{u_k}}$, $u_k\in {\cal U}_{k,d}$, $1\leq k\leq s$,
 of $\sigma$-smooth functions $f_{u_k}(\tb_{u_k})$, where $\tb_{u_k}=(t_{j_1},\ldots,t_{j_k})\in\mathbb{R}^k$,
 $u_k=\{j_1,\ldots,j_k\}$, $1\leq j_1< \ldots< j_k\leq d$, as defined in Section \ref{RC},
it might be of interest to consider the balls of analytic
 functions $f_{u_k}(\tb_{u_k})$, $\tb_{u_k}=(t_{j_1},\ldots,t_{j_k})\in\mathbb{R}^k$,
 with smoothness parameter $\sigma>0$
such that (a) $f_{u_k}$ is 1-periodic in each of its arguments, (b) $f_{u_k}$ can be analytically continued from $\mathbb{R}^k$ to
the $k$-dimensional strip $\mathbb{S}_\sigma=\{{\bf z}_{u_k}\in \mathbb{C}^k: |{\rm Im}\, {\bf z}_{u_k}|\leq \sigma/(2\pi)\}$, and $|f({\bf z}_{u_k})|\leq M$
for all ${\bf z}_{u_k}\in \mathbb{S}_{\sigma}$
and some $M>0$. A more general setup, when the regularity constraints are imposed on the whole $d$-variate function $f$ rather than
on its $k$-variate components  for $1\leq k\leq s$,
may also be addressed.

Variable selection techniques for high-dimensional data are widely applied across numerous research fields. In the context of functional ANOVA model (in practice also called Smoothing Spline ANOVA decomposition, or SS-ANOVA), the problem of variable selection was studied, for example, in Wahba et al. (1995) with applications to data from health sciences. Therefore, it is also of interest to investigate how the selection method proposed in this article may be applied to health and well-being data.

\section*{Acknowledgements} The research of N. Stepanova was supported by a Discovery Grant of the
Natural Sciences and Engineering Research Council of Canada. The research of M. Turcicova was partially supported by Discovery Grants of the  Natural Sciences and Engineering Research Council of Canada during the author's stays at Carleton University in the 2023--2024 academic year. Additionally, the author's research was funded by the project ``Research of Excellence on Digital Technologies and Wellbeing CZ.02.01.01/00/22\_008/0004583" which is co-financed by the European Union.

\section*{Disclosure statement}
The authors confirm that there are no relevant financial or non-financial competing interests to report.

\section*{Data declaration}
This paper does not report any research data or analysis.

\end{document}